\documentclass[11pt,a4paper]{amsart} 
\usepackage[english]{babel}
\usepackage{amssymb,xspace,enumerate,graphics,ifthen}
\usepackage{amstext}
\theoremstyle{plain}
\usepackage{amsbsy,amsfonts,latexsym,enumerate}

\newtheorem{theorem}{Theorem}[section]
\newtheorem{lemma}[theorem]{Lemma}
\newtheorem{proposition}[theorem]{Proposition}
\newtheorem{corollary}[theorem]{Corollary}

\newtheorem{definition}[theorem]{Definition}

\def \pint {\vbox{ \hbox to 5 pt {\hfil \vrule height 4pt}\hrule}\hskip 3pt}

\def\leqs{\lesssim}

\def\geqs{\gtrsim}
\def\eqs{\eqsim}
\def\cc{\mathbb{C}}
\def\rr{\mathbb{R}}

\def\cb{{\mathcal{B}}}
\def\ch{{\mathcal{H}}}
\def\cn{{\mathcal {N}}}
\def\cp{{\mathcal {P}}}

\def\cv{{\mathcal {V}}}
\def\co{{\mathcal{O}}}

\def\supp{{\rm supp}}
\def\vol{{\rm Vol}}
\def\sp{{\rm sp}}

\def \qed {\hbox{\hskip 5pt} \vbox{\hrule \hbox to 5pt 
{\vrule height 4.2pt \hfil \vrule}\hrule}}

\newcommand{\hl}{{h_\Lambda(z,t)}}
\newcommand{\deltar}{d}
\newcommand{\deltap}{d\hskip-3pt\cdot\hskip-3pt}

\newcommand{\dszp}{d(sz')}
\newcommand{\drz}{{d(rz')}}
\newcommand{\dtz}{{d(tz)}}
\newcommand{\dtsz}{{d(tsz')}}
\newcommand{\dz}{{d(z)}}
\newcommand{\dzeta}{d(\zeta)}
\newcommand{\dzetaO}{d(\zeta_0)}
\newcommand{\dhl}{d(h_\Lambda(z,t))}
\newcommand{\dxi}{d(\xi)}
\renewcommand{\i}{(I)} \newcommand{\ijk}{(I)_{j,k}}
\newcommand{\ii}{(II)}
\newcommand{\iijk}{(II)_{j,k}}
\newcommand{\iijkd}{(II)_{j,k,\delta}}
\newcommand{\iijke}{(II)_{j,k,\varepsilon}}
\newcommand{\iii}{(III)}
\newcommand{\iiijk}{(III)_{j,k}}
\def \pint {\vbox{ \hbox to 5 pt {\hfil \vrule height 4pt}\hrule}\hskip 3pt}
\newcommand{\wrt}{with respect to }

\newcommand{\mat}[3]{\left #1  #2  \right #3}
\newcommand{\diffp}[2]{\frac{\partial #1}{\partial #2}}
\newcommand{\mlabel}[1]{\label {#1}}
\renewcommand{\over}[2]{\genfrac{}{}{0pt}{}{#1}{#2}}

\title{Zero sets of $\ch^p$ functions in convex domains of finite type}
\author{William ALEXANDRE}
\address{Laboratoire Paul Painlev\'e U.M.R. CNRS 8524, U.F.R. de
Math\'ematiques,  cit\'e scientifique, Universit\'e Lille 1, F59 655 Villeneuve d'Ascq Cedex, France.}
\email{ william.alexandre@math.univ-lille1.fr}

\begin{document}

\begin{abstract}
We give a condition under which a divisor $\hat X$ in a bounded convex domain of finite type $D$ in  $\cc^n$ is the zero set of a function in a Hardy space $\ch^p(D)$ for some $p>0$. This generalizes Varopoulos' result
[Zero sets of $\ch^p$ functions in several complex variables, Pac. J. Math. (1980)]
on zero sets of $\ch^p$-functions in strictly convex domains of $\cc^n$.
\end{abstract}
\maketitle
\section{Introduction and main result}
We present here our main result and the outline of its proof.
\subsection{Zero sets of functions in the Nevanlinna and Hardy classes}
We denote by $D=\{z\in\cc^n,\ r(z)<0\}$ a bounded domain in $\cc^n$, where $n$ is a positive integer and $r$ is a smooth function such that $d r\neq 0$ on the boundary of $D$. We set $d=|r|$ and $D_\varepsilon =\{z\in\cc^n,\ r(z)<\varepsilon \}$. We denote by $bD_\varepsilon$ the boundary of $D_ \varepsilon$, by $T^\cc_z bD_{r(z)}$  the complex tangent space to $bD_{r(z)}$ at $z$, and  by $d\sigma_\varepsilon $ the euclidean area measure on $bD_\varepsilon $. The Nevanlinna class $\cn(D)$ is the set of holomorphic functions $f$ on $D$ such that 
$$\sup_{\varepsilon >0}\int_{bD_{-\varepsilon} }\bigl|\log|f(z)|\bigr|d\sigma_{-\varepsilon} (z)<+\infty.$$
The Hardy space $\ch^p(D)$, $p>0$, is the set of holomorphic functions $f$ on $D$ such that
$$\|f\|_p=\left(\sup_{\varepsilon >0}\int_{bD_{-\varepsilon} }|f(z)|^pd\sigma_{-\varepsilon} (z)\right)^{\frac1p}<+\infty.$$

Let $X=\{z\in D,\ f(z)=0\}$ be the zero set of a function $f\in \cn (D)$, let $X_k$ be the irreducible components of $X$ and let $n_k$ be the corresponding multiplicities of $f$ ; the data $\hat X=\{X_k,n_k\}_k$ is commonly called a divisor. It is well known that $X$ or equivalently $\hat X$ satisfies the {\em Blaschke Condition} :
\begin{align}
 \sum_{k} n_k\int_{X_k} \dz d\mu _{X_k}(z)<+\infty. \tag{B}\mlabel{B}
\end{align}
When $D$ is the unit disk of $\cc$, this condition simply becomes the well known condition $\sum_k 1-|a_k|<+\infty$ where $X=\{a_k\}_k$, each $a_k$ counted accordingly to its multiplicity. It is also well known that any sequence $(a_k)_k$ satisfying the Blaschke Condition is the zero set of a function $f$ belonging to $\cn(D)$ and of a function $g$ belonging to $\ch^p(D)$, $p>0$. This in particular means that the functions of the Nevanlinna class and the functions of the Hardy spaces have the same zero sets.

In $\cc^n$, $n>1$, the situation is much more intricate. It was proved independently  by Henkin \cite{Hen} and Skoda \cite{Sko} that when $D$ is strictly pseudoconvex and satisfies some obvious topological condition, any divisor which satisfies the Blaschke Condition (\ref{B}) is the zero set of a function $f\in\cn(D)$. Some partial results are known for the polydisc (\cite{And,BC}), or special domains (\cite{Cha}) and the Henkin-Skoda Theorem was also proved for pseudoconvex domains of finite type in $\cc^2$ (\cite{CNS}), for convex domains of finite strict type in $\cc^n$ (\cite{BCD}) and for convex domain of finite type in $\cc^n$ (\cite{DM}).

In the case of Hardy spaces in $\cc^n$, $n>1$, the situation is even more complicated. Contrary to the one dimensional case, the zero sets of functions in the Nevanlinna class and the zero sets of functions in the Hardy classes are different. Moreover, for distinct $p$ and $q$, the zero sets of functions of $\ch^p$- and $\ch^q$-classes are different (see \cite{Rud}). However, Varopoulos managed to give in \cite{Var} a general condition for a divisor $\hat X$ to be the zero set of an holomorphic function belonging to $\ch ^p(D),$ for some $p>0$. Varopoulos' proof was simplified  by Andersson and Carlsson in \cite{AC}. Bruna and Grellier attempted in \cite{BG} to generalize Varopoulos result to the case of convex domains of finite strict type, but there are some gaps in their proof. We aim to prove in this article the generalization of Varopoulos result to the case of convex domains of finite type in $\cc^n$, which includes in particular the case of convex domains of finite strict type.

\subsection{Varopoulos' result}
We will now present Varopoulos' result and the scheme of its proof that we translate to the framework of convex domains of finite type. We will also explain the differences with the situation of convex domains of finite type.
\par\medskip
Varopoulos used the Lelong current associated with a divisor $\hat X$ in order to define what he called a Uniform Blaschke Condition. Lelong proved that any divisor $\hat X$ can be associated with a closed positive $(1,1)$-current $\theta = \theta_{\hat X}$ of order 0, that is a $(1,1)$-form $\theta=\sum_{j,k=1}^n\theta  _{j,k} dz_j\wedge d\overline{z_k}$, where each $\theta  _{j,k}$ is a complex measure such that $d\theta  =0$ and for all $\lambda _1,\ldots,\lambda _n\in\cc$, $\sum_{j,k=1}^n\theta  _{j,k}\lambda _j\overline\lambda _k$ is a positive measure. 
The Blaschke Condition $(\ref{B})$ can be reformulated by asking that $\deltap  |\theta  |$ be a bounded measure on $D$. Varopoulos condition also involved $\deltap\theta$, and in particular required $\deltap |\theta  |$, to be not only a bounded measure, but a Carleson measure. 
We here give the definition of Carleson measures in the setting of convex domains of finite type. This notion is related to the structure of homogeneous space on $D$ induced by the polydics of McNeal defined in \cite{McN, McN1, McN2}. They are the analog of Koranyi balls of strictly convex domains and they are defined as follows. For $z$ near $bD$, small positive $\varepsilon $ and $v\in\cc^n$, $v\neq 0$, we set
\begin{align*}
\tau(z,v,\varepsilon)&:=\sup\{t>0, |r(z+\lambda v)-r(z)|<\varepsilon, \forall \lambda \in\cc,\ |\lambda |<t\}.
\end{align*}
This positive number $\tau(z,v,\varepsilon)$ is the distance from $z$ to the level set $\{r=r(z)+\varepsilon\}$ in the complex direction $v$.
We now recall the definition of an $\varepsilon$-extremal basis $w_1^*,\ldots, w_n^*$ at the point $z$, given in \cite{BCD} :  $w_1^* =\eta_z$ is the outer unit normal to $bD_{r(z)}$ at $z$ and if $w^*_1,\ldots, w^*_{i-1}$ are already defined, then $w_i^*$ is a unit vector orthogonal to $w_1^*,\ldots, w^*_{i-1}$ such that $\tau(z,w^*_i,\varepsilon)= \sup_{\over{v\perp w^*_1,\ldots, w^*_{i-1}}{\|v\|=1}} \tau(z,v,\varepsilon)$. When $D$ is strictly convex, $w_1^*$ is the outer unit normal to $bD_{r(z)}$ and we may choose any basis of $T^\cc_z bD_{r(z)}$ for  $w^*_2,\ldots, w^*_n$. Therefore, when $D$ is strictly convex, an $\varepsilon$-extremal basis at $z$ can be chosen smoothly depending on the point $z$. Unfortunately, this is not the case for convex domains of finite type (see \cite{Hef2}).\\
We put $\tau_i(z,\varepsilon)=\tau(z,w_i^*,\varepsilon)$, for $i=1,\ldots,n$.
Writing $A\leqs B$ if there exists a constant $c>0$ such that $A\leq cB$ and $A\eqs B$ if $A\leqs B$ and $B\leqs A$ both hold, we have for a strictly convex domain $\tau_1(z,\varepsilon)\eqs \varepsilon$ and $\tau_j(z,\varepsilon)\eqs\varepsilon^{\frac12}$ for $j=2,\ldots,n$. For a convex domain of finite type $m$, we only have $\varepsilon^{\frac12}\leqs \tau_n(z,\varepsilon)\leq\ldots\leq \tau_n(z,\varepsilon)\leqs \varepsilon^{\frac 1m}$,  uniformly with respect to $z$ and $\varepsilon$.\\
The McNeal polydisc centered at $z$ of radius $\varepsilon $ is the set $$\cp_\varepsilon(z):=\left\{ \zeta=z+\sum_{i=1}^n \zeta^*_iw_i^*\in\cc^n, |\zeta^*_i|<\tau_i(z,\varepsilon),\:i=1,\ldots, n\right\}.$$ 
\begin{definition}
We say that a positive finite measure $\mu$ on $D$ is a {\it Carleson measure} and we write $\mu\in W^1(D)$ if 
$$\|\mu\|_{W^1(D)}:= \sup_{\over{z\in bD}{\varepsilon>0}}\frac {\mu(\cp_\varepsilon (z)\cap D) }{\sigma(\cp_\varepsilon (z)\cap bD )}<\infty.$$
\end{definition}
Varopoulos Uniform Blaschke Condition requires that
\begin{align}
 \deltap |\theta  |, \quad \deltar ^{\frac12} |\partial r\ \wedge \theta  |, \quad  \deltar ^{\frac12} |\partial r \wedge \theta  |, \quad \text{and} \quad |\partial r\wedge \overline\partial r \wedge \theta  | \text{ belong to } W^1(D). \tag{UB}\mlabel{UB}
 \end{align}
The factors $\deltar $, $\deltar ^{\frac12}$ are weights which actually depend on the components of $\theta  $. For example in $\partial r\wedge \theta  $, the exterior product of $\theta  $ with $\partial r$ cancels the normal component of $\theta  $ in $dz$ so that only the tangential part of $\theta  $ in $d\overline z$ is left. Varopoulos put in front of this tangential part a factor $\deltar ^{\frac12}$, the exponent $\frac12$ being, in Varopoulos' case  of strictly pseudoconvex domains, $1$ over the order of contact of a tangent vectors field and the boundary of $D$. This in particular means that the normal component of $\theta  $ can behave in a  worse manner than the $dz$-tangential component which itself can behave in a  worse manner than the whole tangential component, this worse behavior being quantified by the order of contact of vectors fields with the boundary.

The situation is more complicated in the case of convex domains of finite type because the order of contact of tangential vectors fields is not constant. In order to overcome this difficulty, we use the following norm defined in \cite{BCD}.
For $z\in\cc^n$ and $v$ a non zero vector we 
set
\begin{align*}
k(z,v)&:=\frac{\dz }{\tau(z,v,\dz )}.
\end{align*}
For a fixed $z$, the convexity of $D$ implies that the function defined by $v\mapsto k(z,v)$ if $v \neq 0$, $0$ otherwise, is a kind of non-isotropic norm. In the case of strictly convex domains, when $v$ belongs to $T^\cc_z bD_{r (z)}$, $k(z,v)$ is comparable to $\deltar ^{\frac12}$, whereas if $\eta_z$ is the unit outer normal to $bD_{\dz)}$ at $z$, $k(z,\eta_z)$ is  comparable to $1$. This implies that the factor $\deltar ^{\frac12}$ in (\ref{UB}) is equal to $\frac{\deltar }{k(\cdot,v)k(\cdot,\eta_z)}$ where $v$ is any tangent vector field, the factor $\deltar $ is in fact $\frac{\deltar }{k(\cdot,\eta_z)k(\cdot,\eta_z)}$ and the factor ``1'' in front of $|\partial r\wedge\overline{\partial} r\wedge \theta|$ is actually $\frac{\deltar }{k(\cdot,v)k(\cdot,w)}$, $v$ and $w$ being any tangent vector fields.

In a Uniform Blaschke Condition for convex domains of finite type, $v$ and $w$ have to appear explicitly because we need to link the weight $\frac{\deltar }{k(\cdot,v)k(\cdot,w)}$ and the ``component of $\theta  $ in the directions $v$ and $w$''.
\subsection{Main result}
In order to have a Uniform Blaschke Condition type which makes sense for general currents with measure coefficients and not only for smooth currents, we set the following definition (compare with the Uniform Blaschke Condition of \cite{BG} which makes sense only for smooth currents) :
\begin{definition}
We say that a $(p,q)$-current $\mu$ of order $0$ with measure coefficients  is a $(p,q)$-Carleson current if
$$\|\mu\|_{W^1_{p,q}}:=\sup_{u_1,\ldots,u_{p+q}} \left\|\frac1{k(\cdot,u_1)\ldots k(\cdot,u_{p+q})}\left|\mu(\cdot)[u_1,\ldots,u_{p+q}]\right|\right\|_{W^1}<\infty,$$
where the supremum is taken over all smooth vector fields $u_1,\ldots, u_{p+q}$ which never vanish and where $|\mu(\cdot)[u_1,\ldots,u_{p+q}]|$ is the absolute value of the measure $\mu(\cdot)[u_1,\ldots,u_q]$.\\
We denote by $W^1_{p,q}(D)$ the set of all $(p,q)$-Carleson currents.\\
A $r$-Carleson current is a sum of $(p,q)$-Carleson currents with $p+q=r$.
\end{definition}
This norm was already defined and used in \cite{WA1}. It is a norm on forms with measure coefficients associated with the vectorial norm $k$. It is defined in the same spirit as the norms used in \cite{AB} and \cite{AC} but $\|\cdot\|_{W^1_{p,q}}$ takes into account the non isotropy of the boundary of the domain. We should also notice that our norm is weaker than the norm of Bruna, Charpentier and Dupain (see \cite{BCD}) in the sense that up to a uniform multiplicative constant, $\|\mu\|_{W^1_{(1,1)}}$ is bounded by $ \int_B\|\mu(\zeta )\|_kd\lambda(\zeta)$ for all smooth $\mu$, where $\|\mu(\zeta)\|_k:=\sup\left\{\frac{|\mu(u)|}{k(\zeta,u)},\ u\neq 0\right\}$. Moreover, we point out that $\|\cdot\|_{W^1_{(p,q)}}$ is defined for smooth currents but also for general currents with measure coefficients.
Finally, in the case of strictly convex domains, we notice that if $\theta$ satisfies the Uniform Blaschke Condition (\ref{UB}), then $\theta  $ is a $(1,1)$-Carleson current.
We will prove the following theorem which is our main result and which, together with  the preceding remarks, extends Varopoulos' result \cite{Var} to the case of convex domains of finite type :
\begin{theorem}[{\bf Main Theorem}]
 Let $D$ be a $C^\infty$-smooth convex domain of finite type, $\hat X$ a divisor in $D$, $\theta  _{\hat X}$ the $(1,1)$-current of Lelong associated with $\hat X$. Then, if $d\cdot \theta_{\hat X}$ is a Carleson current, there exist $p>0$ and $f\in\ch^p(D)$ such that $\hat X$ is the zero set of $f$.
\end{theorem}

\subsection{Scheme of the proof of the main result}
The main scheme of the proof is classical~: we have to find a real valued function $u$ such that $i\partial\overline\partial u = \theta  _{\hat X}$ with a growth condition on $u$. Since $D$ is convex, such a function $u$ is equal to $\log|f|$ for an $f$ that defines $\hat X$.\\
In order to find $u$, we proceed in two steps. First we solve the equation $idw=\theta  _{\hat X}$ with $w$ such that $w=-\overline w$. This is done thanks to the following theorem.

\begin{theorem}\mlabel{th1}
 Let $D$ be a $C^\infty$-smooth bounded convex domain of finite type, $\theta$ a $d$-closed $(1,1)$-current  of order $0$ such that $\deltap \theta$ is a Carleson current. Then there exists a real $1$-Carleson current $\omega$ such that $d\omega=\theta  $.
\end{theorem}

We then set $w=-i\omega$, where $\omega$ is given by Theorem \ref{th1}, so that $\overline{w}=-w$.  We write $w$ as $w=w_{1,0}+w_{0,1}$ where $w_{1,0}$ is a $(1,0)$-Carleson current and $w_{0,1}$ is a $(0,1)$-Carleson current. We trivially have $\overline w_{0,1} = -w_{1,0}$. Moreover, since $idw= \partial w_{1,0} + \overline \partial w_{1,0} + \partial w_{0,1} +\overline\partial w_{0,1}$, and since $idw=\theta  $ is a $(1,1)$-current, for bidegree reasons we have $\partial w_{1,0}=0$, $\overline \partial w_{1,0}= \partial w_{0,1}$ and $\overline\partial w_{0,1}=0$. Since $D$ is convex, we can find $v$ such that $\overline{\partial} v=w_{0,1}$. Setting $u=2\Re v$, we get
\begin{align*}
 i\partial\overline\partial u&= i\partial\overline\partial v- i\overline\partial\partial \overline v\\
 &= i\partial\overline\partial v- i\overline\partial\ \overline{\overline\partial v}\\
 &= i\partial w_{0,1} -\overline \partial \overline{w_{0,1}}\\
 &= i\partial w_{0,1} +\overline \partial w_{1,0}\\
 &= idw =\theta.
\end{align*}
Therefore, in order to prove our main theorem, we have to find a solution of the $\overline\partial$-equation $\overline\partial v=w_{0,1}$ with $\exp v$ in $L^p(bD)$. It is given by the following theorem :
\begin{theorem}\mlabel{th2}
 Let $D$ be a $C^\infty$-smooth convex domain of finite type and let $\omega$ be a $\overline\partial$-closed $(0,1)$-Carleson current in $D$. Then there exist $p>0$ and a solution $v$ to the equation $\overline\partial v=\omega$ such that $\exp v$ belongs to $L^p(bD)$.
 \end{theorem}
 
 We now give the scheme of the proofs of Theorems \ref{th1} and \ref{th2}. In order to prove Theorem \ref{th1}, without restriction, we will assume that $0$ belongs to $D$, that $0$ does not belong to $\supp(\theta)$ and that $\theta$ is supported in a sufficiently small neighborhood of $bD$. We will use the Poincar\'e homotopy operator and we need a deformation retract $h:D\times[0,1]\to D$ of $D$ onto 0.
 Using convexity, Bruna, Charpentier and Dupain simply defined $h$  by $h(z,t)=t\cdot z$. However, as already pointed out by Varopoulos in \cite{Var}, this choice does not work for Hardy spaces. In this case, it is necessary to take the mean value of a suitable family of homotopy operators. We now give an analogue for strictly convex domains of the deformation retract used by Andersson-Carlsson \cite{AC} and Varopoulos \cite{Var}.
 \par\medskip
Still assuming that $0$ belongs to $D$, we denote by $p$ the calibrator or gauge function for $D$, that is $p(\zeta  )=\inf\{\lambda >0, \  z\in \lambda D\},$ and from now on $r=p-1$. We notice that since $p$ is homogeneous, the level sets $bD_{\varepsilon}$ are homotetic and $\tau (z,v,\varepsilon )$ itself becomes homogeneous. Moreover, with such a choice of a defining function, for all $t>0$, any $v$ belongs to $T^\cc_zbD_{r(z)}$ if and only if it belongs to $T^\cc_{tz}bD_{r(tz)}$.
 
Let $w_1^*(z)$ be the outer unit normal to $bD_{r(z)}$ at $z$, let $w_2^*(z),\ldots, w_n^*(z)$ be a basis of $T^\cc_zbD_{r(z)}$ smoothly depending on $z$, which is always possible at least locally. We notice that, for all $t>0$, $w_1^*(z)$ is the outer unit normal to $bD_{r(tz)}$ at $tz$, and that $w_2^*(z),\ldots, w_n^*(z)$ is a basis of $T^\cc_{tz}bD_{r(tz)}$. Therefore we can assume that $w_j^*(tz)=w_j^*(z)$ for all $t>0$. Then for $\Lambda=(\lambda _1,\ldots,\lambda _n) \in \Delta^n$, $\Delta=\{\xi\in\cc,\ |\xi|<1\}$, define $h_\Lambda :D\times[0,1]\to D$ by
$$h_\Lambda(z,t)=tz+t\left((1-t) \lambda_1 w^*_1(tz)+ \sum_{j=2}^n \lambda_j \sqrt{1-t}\cdot w_j^*(tz)\right),$$
and set
$$H\theta  =\frac1{(\vol(\Delta))^n}\int_{\Delta^n} \left(\int_{[0,1]} h^*_\Lambda \theta  \right)d\Lambda $$
where the inner integral is the $t$-integral of the $dt$-component of $h^*_\Lambda \theta  $.
We have $h_\Lambda(z,0)=0$ and $h_\Lambda(z,1)=z$ for all $z\in D$, $h_\Lambda$ is smooth in $D\times ]0,1[$ for all $\Lambda$ and thus $dH+Hd=Id$. Let us look a bit at what $h_\Lambda $ and $H$ do. When $\Lambda $ is fixed, $h_\Lambda (z,\cdot)$ is a path from $0$ to $z$ which, for all $\Lambda\neq 0$, is not a straight line as in \cite{BCD}. Each $h_\Lambda $ induces an homotopy operator and $H$ is in fact the mean value of these homotopy operators.\\
Let us fix $z$ and $t$ and let $\Lambda $ varies over $\Delta^n$. When $D$ is a strictly convex domain, the factor $\sqrt{1-t}$ in $h_\Lambda$ is comparable to $\tau(tz,w_j^*(tz),1-t)$ for all $j=2,\ldots, n$, the factor $1-t$ is comparable to  $\tau(tz,w_j^*(tz),1-t)$. In particular, when $\Lambda $ varies over $\Delta ^n$, the image of $h_\Lambda (z,t)$ is $\cp_{1-t}(tz)$.

So, when $D$ is a convex domain of finite type, our first attempt at a proof could simply be to replace in $h_\Lambda $ the vectors $w^*_j(tz)$ by a $(1-t)$-extremal basis at $tz$ that we still denote by $w^*_j(tz)$, $j=1,\ldots, n$,  and the factor $\sqrt{1-t}$ by $\tau (tz,w^*_j(tz),1-t)$ for $j=2,\ldots, n$, $1-t$ by $\tau (tz,w^*_1(tz), 1-t)$. However, $h_\Lambda $ would not be smooth because $\varepsilon$-extremal bases at $z$ may behave in a really bad way and in general, do not depend continuously on $\varepsilon$ or on $z$ (see \cite{Hef2}). We have to find a smooth way of describing $\cp_{1-t}(tz)$.
More precisely, we look for a smooth map $h_\Lambda:D\times[0,1]\to D$ with the following properties : for all $\Lambda$ in $\Delta_n(\rho)=\{\Lambda\in\cc^n,\ |\Lambda|<\rho\}$ (where $\rho>0$ is a small number which has  to be determined) $h_\Lambda(z,0)=0$, $h_\Lambda(z,1)=z$, and there exist a uniform constant $\gamma>0$, and $C>c>0$ depending on $\rho$ such that for fixed $z\in D$, $t\in[0,1-\gamma\dz ]$ :
$$c\cp_{1-t}(tz)\subset \{h_\Lambda(z,t),\ |\Lambda|<\rho\}\subset C\cp_{1-t}(tz).$$
We will explain later why we  only require that these properties hold only for $t\in [0,1-\gamma\dz ]$ and not for $t$ in the whole interval $[0,1]$.
Moreover, for technical reasons that will become clear later on, we also want that $C$ goes to $0$ when $\rho$ goes to $0$.\\
We will achieve this aim thanks to the Bergman metric (see Subsection \ref{ssI.1} for the definition of the Bergman metric). The next two propositions link McNeal polydiscs and the Bergman metric in convex domains of finite type. The first one  was proved by McNeal in \cite{McN1}.
\begin{proposition}\mlabel{prop0.1}
 Let $\zeta\in D$ be a point near $bD$, $\varepsilon>0$ and $w^*_1,\ldots, w^*_n$ an $\varepsilon$-extremal basis at $\zeta$ and $v= \sum_{j=1}^n v^*_j w^*_j$ a unit vector. Then, uniformly \wrt $\zeta  , v$ and $\varepsilon $, we have
\begin{align*}
\frac1{\tau (\zeta  ,v,\varepsilon )}&\eqs  \sum_{j=1}^n\frac{|v^*_j|}{\tau_j(\zeta  ,\varepsilon )}.
\end{align*}
\end{proposition}
Therefore $\cp_\varepsilon (\zeta  )$ could also be defined as the set $\{\zeta  +\lambda v,\ v\in\cc^n, |v|=1, \lambda\in\cc, |\lambda|<\tau(\zeta  ,v,\varepsilon )\}$.\\
Now, let $B(\zeta  )$ be the matrix in the canonical basis which determines the Bergman metric $\|\cdot\|_{B,\zeta  }$ at $\zeta  $, i.e. $\|v\|_{B,\zeta  }=\overline v^t B(\zeta  ) v$ for any vector $v$. We recall that $B$ depends  smoothly on $\zeta\in D$ but explodes on the boundary. The following result was proved by McNeal in \cite{McN2}.
\begin{proposition}\mlabel{propI.1.7}
Let $\zeta  \in D$ be a point near $bD$, $v$ a unit vector in $\cc^n$. Then, uniformly with respect to $\zeta  $ and $v$,
$$\|v\|_{B,\zeta  }\eqs \frac1{\tau (\zeta  ,v,\dzeta)}.$$
\end{proposition}
Therefore there exist $C>c>0$ such that for all $\zeta  $ near $bD$ and $\rho>0$
\begin{align}
c\rho  \cp_{\dzeta}(\zeta  )\subset \{\zeta  +\lambda v \in\cc^n,\ |v|=1 \text{ and }  \|\lambda v\|_{B,\zeta  }<\rho\} \subset C\rho \cp_{\dzeta}(\zeta  ).\mlabel{eqi1} 
\end{align}

Since the Bergman metric is an hermitian metric, for all $\zeta  \in D$ there exists a positive hermitian matrix $A(\zeta  )$ such that $A(\zeta  )^{-2}=B(\zeta  )$. The inverse mapping theorem ensures that $A$ depends smoothly on $\zeta  \in D$, and $\|A(\zeta  )v\|_{B,\zeta   }=|v|$ for all $\zeta  $ and $v$.
Therefore $(\ref{eqi1})$ becomes
\begin{align}
c\rho  \cp_{\dzeta}(\zeta  )\subset \{\zeta  +A(\zeta  ) v,  |v|<\rho\} \subset C\rho \cp_{\dzeta}(\zeta  ).\mlabel{eqi1ter}
\end{align}
Putting $\zeta=tz$, since $\dtz=1-t+t\dz \eqs 1-t$ when $t\leq 1-\gamma\dz$, Corollary \ref{corI.1.2} yields
\begin{align}
c\rho  \cp_{1-t}(tz)\subset \{tz+A(tz) v,  |v|<\rho\} \subset C\rho \cp_{1-t}(tz)\mlabel{eqi1bis} 
\end{align}
for all $z$ and $t$ such that $0\leq t\leq 1-\gamma \dz$. In other words,  $\{tz+A(tz) v,  |v|<\rho\}$ is almost equal to 
$\rho \cp_{1-t}(tz)$.

For $t$ close to 1, we cannot use $A$ in order to get a set which is almost equal to $\cp_{1-t}(tz)$. Indeed, $A(\zeta)$ yields a set which is almost equal to $\cp_{\dzeta}(\zeta  )$ and by homogeneity of $D$, it is possible to obtain a set which is almost equal to $  \cp_{1-t}(tz)$ using $A(\zeta  )$ as in (\ref{eqi1bis}) with a point $\zeta=\lambda z  \in D$ such that $\dzeta\eqs1-t$ (that is for a point $\zeta  $ close to $bD$ if $t$ is close to $1$). However, when $\zeta  $ goes to the boundary, the derivatives of $A(\zeta  )$ explode, and actually they explode so much and the computations will not work. This problem does not appear in the strictly convex case because the extremal bases can be chosen to be smooth in a neighborhood of $\overline D$.

It appears in the computations that, when $1-t\leq \gamma \dz$, there is in fact no need to take mean value of homotopy operators. But, in order that things work when $1-t\eqs \dz$, we have to make a cleverer choice of retracts. We define $h_\Lambda$ as follows. Let $\varphi $ be a $C^\infty$ smooth function such that $\varphi (t)=1$ if $t<\frac12$, $\varphi (t)=0$ if $t>1$, and define the map $h_\Lambda :D\times [0,1]\to D$ for $|\Lambda|\leq \rho$ by
\begin{align*}
h_\Lambda (z,t)
&=tz+ t\varphi \left(\frac{1-t}{\gamma\dz}\right)\frac{1-t}{\dz} A(z) \cdot\Lambda   +t\left(1-\varphi \left(\frac{1-t}{\gamma\dz}\right)\right)A(tz) \cdot\Lambda 
\end{align*}
where $\gamma$ has to be chosen sufficiently small.\\
The associated homotopy operator is
$$H\theta  =\frac1{\vol(\ \Delta_n(\rho))}\int_{\Lambda \in \Delta_n(\rho)} \left(\int_{t\in [0,1]} h^*_\Lambda \theta  \right)d\Lambda. $$
The map $h_\Lambda $ is $C^\infty$-smooth in $D\times ]0,1[$, $h_\Lambda (z,0)=0$ and $h_\Lambda (z,1)=z$ for all $z$ in $D$.
For fixed $z$ and $t$ such that $1-t\geq \gamma\dz$ we get from (\ref{eqi1bis})
\begin{align}
ct\rho  \cp_{1-t}(z)\subset \{h_\Lambda(z,t),\ |\Lambda|<\rho\}\subset Ct\rho \cp_{1-t}(z).\mlabel{1-t>delta}
\end{align}
From (\ref{eqi1ter}), for fixed $z$ and $t$ such that $1-t\leq \gamma\dz$ we have
\begin{align}
\{h_\Lambda(z,t),\ |\Lambda|<\rho\}\subset C\rho\cp_{\dz }(z).\mlabel{1-t<delta}
\end{align}

Now that we have obtained a good homotopy formula, the rest of the proof of Theorem \ref{th1} consists of tedious computations that we carry out in Section \ref{section2}. In order to estimate $H\theta  $, we will distinguish three cases, depending on whether $1-t\leq \frac \gamma2\dz$, $1-t\geq \gamma\dz$ or $\frac \gamma2\dz\leq 1-t\leq \gamma\dz$ (see Subsections \ref{ssI.3}, \ref{ssI.4} and \ref{ssI.5} respectively).\\
We will be led to compute derivatives of $h_\Lambda$ and so of $A$. We will compute these derivatives by applying  the inverse mapping theorem to the map $\Phi$ defined on the set of positive hermitian matrices by $\Phi(B)=B^{-2}$. In order to compute $d\Phi^{-1}$, we will have to solve the equation $BM+MB=M'$ where $M'$ is given and where $M$ is an unknown matrix. Because we need optimal estimates, we will give an explicit expression of $M$ using ideas of Rosenblum \cite{Ros}. This will be done in Subsection \ref{ssI.2}, after we have given in Subsection \ref{ssI.1} the tools related to convex domains of finite type.
\par\medskip
The proof of Theorem \ref{th2} is more classical. We will follow ideas of \cite{AC}, \cite{Sko} and \cite{DM} that we have to adapt to our new norm $\|\cdot\|_{W^1}$. We will use Diederich-Mazzilli's solution of the $\overline\partial$-equation, which itself involved a Skoda type integral operator constructed with the Dierderich-Forn\ae ss support function $S$ for convex domains of finite type. In order to prove Theorem \ref{th2}, we will have to estimate the $W^1$-norm of our solution. Therefore we will need to find suitable vectors fields. It turns out that extremal bases realized the supremum in the kind of norm $\|\cdot\|_k$ used in \cite{BCD}. However, we need here smooth vectors fields and as we already said, extremal bases are not smooth. The Bergman metric (again) will give us vectors fields which will be a smooth alternative to extremal bases (see Section \ref{section3} for details).

\section{The $d$-equation}\mlabel{section2}
In order to prove Theorem \ref{th1}, we have to prove that for all non-vanishing vector fields $u$, all $z_0\in bD$, all $\varepsilon >0$, the following inequality holds uniformly :
\begin{align}
 \int_{\cp_\varepsilon(z_0) \cap D} \frac1{k(z,u(z))} |H\theta  (z)[u(z)]| d\lambda(z) &\leqs \sigma(\cp_\varepsilon (z_0)\cap bD)\|\theta  \|_{W^1_{1,1}}.\mlabel{eqI.5}
\end{align}
By standard regularization arguments (see \cite{AC}), we can assume $\theta$ smooth on $D$.  
When we compute $H\theta  (z)[u(z)]$, we get
\begin{align}
\mlabel{eqI.5bis} H\theta  (z)[u(z)]&=\frac1{\vol(\Delta_n(\rho))}\int_{\Lambda \in \Delta_n(\rho)}\int_{t\in[0,1]} \theta  (h_\Lambda(z,t))
 \left[\diffp{h_\Lambda}t (z,t),d_zh_\Lambda(z,t)[u]\right] dt d\Lambda.
\end{align}
The definitions of $\|\theta  \|_{W^1_{1,1}}$ and $h_\Lambda $ naturally lead us to compute $k(\zeta  ,dA_\zeta  [u]\cdot\Lambda)$. We will do this in Subsection \ref{ssI.2} after having recalled in Subsection  \ref{ssI.1} 
the tools for convex domains of finite type that we will need in this section. 

As we will see in the next subsection, the properties of convex domain of finite type are known only in a neighborhood of the boundary. This is why, without restriction since $D$ is convex, we assume that $\supp(\theta)\subset \overline{D\setminus D_{-\varepsilon_0}}$, $\varepsilon_0>0$ as small as we want. Moreover, since $|\hl-tz|\leqs \rho$ uniformly with respect to $\rho, t$ and $z$, if $t$ is small enough, $\hl$ does not belong to $\supp(\theta)$. Therefore there exists a uniform $t_0>0$ such that we only integrate in (\ref{eqI.5bis}) for $t\in[t_0,1]$.

\subsection{Some tools for convex domains of finite type}\mlabel{ssI.1}
We collect here many of the properties of McNeal's polydiscs and of the radii $\tau(z,v,\varepsilon)$. The first ones come
directly from their definition :
\begin{proposition}\mlabel{propI.1.6}
 For all $v\in\cc^n$, all $\zeta\in D$, all $\varepsilon >0$ and all $\lambda\in\cc^*$ :
 $\tau(\zeta,v,\varepsilon)= |\lambda|\tau(\zeta,\lambda v,\varepsilon).$\\
 If $v$ is a unit vector belonging to $T^\cc_\zeta bD_{r(\zeta)}$, then $\varepsilon ^{\frac12}\leqs\tau(\zeta,v,\varepsilon)\leqs \varepsilon^{\frac1m}$, uniformly with respect to $\zeta$, $v$ and $\varepsilon$.\\
 If $v=\eta_\zeta$ is the outer unit normal to $bD_{r(\zeta)}$ at $\zeta$, then $\tau(\zeta,\eta_\zeta,\varepsilon)\eqs \varepsilon$.
\end{proposition}

The next property is proved in \cite{BCD}.
\begin{proposition}\mlabel{propI.1.5}
 Let $z\in D$ be a point near $bD$, $v$ a unit vector in $\cc^n$ and $\varepsilon _1\geq \varepsilon _2>0$. Then we have uniformly with respect to $z$, $\varepsilon _1,\varepsilon _2$ and $v$
$$\left(\frac{\varepsilon_1}{ \varepsilon _2}\right)^{\frac1m}\leqs\frac{\tau (z,v,\varepsilon _1)}{\tau (z,v,\varepsilon _2)}\leqs \frac{\varepsilon _1}{\varepsilon _2}.$$
\end{proposition}
As a corollary of Propositions \ref{prop0.1} and \ref{propI.1.5} we have :
\begin{corollary}\mlabel{corI.1.2}
 Let $z\in D$ be a point near $bD$. If $\varepsilon_1, \varepsilon_2, k, K>0$ are such that  $k\varepsilon _1\leq\varepsilon _2\leq K\varepsilon _1$, there are constants $C\geq c>0$, depending only on $k$ and $K$, such that $$c\cp_{\varepsilon_1}(z)\subset \cp_{\varepsilon _2}(z)\subset C\cp_{\varepsilon _1}(z).$$
 \end{corollary}
In particular, for all $c>0$, $\vol({\cp}_{c\varepsilon}(z))\eqs \vol({\cp}_{\varepsilon}(z))$ uniformly with respect to $z$ and $\varepsilon$.
The following proposition, proved in \cite{McN1}, and Corollary \ref{corI.1.2} show that the polydiscs define a structure of homogeneous space on $D$. 
\begin{proposition}\mlabel{propI.1.9}
 There exists $C>0$ such that, for all $\varepsilon>0$ and all $z,\zeta$ in a neighborhood of $bD$, the following holds true: if ${\cp}_\varepsilon(z)\cap {\cp}_\varepsilon(\zeta)\neq \emptyset$ we have ${\cp}_\varepsilon(z)\subset C{\cp}_\varepsilon(\zeta)$. In particular, $\vol( {\cp}_\varepsilon(z))\eqs \vol( {\cp}_\varepsilon(\zeta))$  uniformly \wrt $\zeta,z$ and $\varepsilon$.
\end{proposition}
We set for $\zeta,z$ near $bD$
$$\delta(z,\zeta):=\inf\{\varepsilon >0,\ \zeta\in {\cp}_\varepsilon(z)\}.$$
Corollary \ref{corI.1.2} and Proposition \ref{propI.1.9} show that $\delta$ is a pseudodistance.

The following proposition is established in \cite{McN1}.
\begin{proposition}\mlabel{propI.1.3}
There exists $c>0$ sufficiently small such that for all $z\in D$ near $bD$, all $\zeta\in c\cp_{\dz }(z)$, we have $\dz \eqs\dzeta$, uniformly with respect to $z$ and $\zeta$ .
\end{proposition}

The following proposition, shown in \cite{McN1}, allows us to compare $\tau(z,v,\varepsilon)$ for different points $z$.
\begin{proposition}\mlabel{propI.1.4}
for all $z\in D$ near $bD$, all  unit vector $v$ in $\cc^n$, all $\varepsilon>0$ and all $\zeta\in\cp_\varepsilon (z)$, we have uniformly with respect to $z$, $\zeta$, $\varepsilon $ and $v$
$$\tau(z,v,\varepsilon)\eqs\tau(\zeta,v,\varepsilon).$$
\end{proposition}
As a corollary of Propositions \ref{propI.1.5}, \ref{propI.1.3} and \ref{propI.1.4}, we have
\begin{corollary}\mlabel{corI.1.3}
There exists $c>0$ such that for all $z$ near $bD$, all $\zeta\in c\cp_{\dz }(z)$, all $v\in\cc^n$ :
$$k(\zeta,v)\eqs k(z,v),$$
 uniformly with respect to $z,\zeta$ and $v$.
\end{corollary}

We will also need the following proposition (see \cite{WA,BCD,DFF}):
\begin{proposition}\mlabel{propI.1.8}
 Let $w$ be any orthonormal coordinates system centered at $\zeta$ and let $v_j$
be the unit vector in the $w_j$-direction. For all multiindices $\alpha$ and $\beta $ with $|\alpha+\beta |\geq 1$ and all $z\in \cp_\varepsilon(\zeta)$:
$$\left|\diffp{^{|\alpha|+|\beta|} r}{w^\alpha\partial \overline w^\beta}(z) \right|\leqs \frac{\varepsilon }{\prod_{j=1}^n \tau (\zeta,v_j,\varepsilon)^{\alpha_j+\beta _j}}
$$
uniformly with respect to $z$, $\zeta$ and $\varepsilon$.
\end{proposition}

We now briefly recall the definition of the Bergman metric (see \cite{Ran}) and  its properties on a convex domain of finite type. The orthogonal projection from $L^2(D)$ onto $L^2(D)\cap {\co}(D)$, where $\co(D)$ is the set of holomorphic function on $D$, is called the Bergman projection. We denote it by $\cb$. There exists a unique integral kernel $b$  such that for all $f\in L^2(D)$ :
$${\cb}f(z)=\int_Db(\zeta,z) f(\zeta )d\lambda(\zeta).$$
The kernel $b(\zeta,z)$ is called the Bergman kernel. This kernel is holomorphic with respect to $z$, antiholomorphic with respect to $\zeta$ and satisfies $b(\zeta ,z)=\overline{b(z,\zeta )}$.\\
The Bergman metric $\|\cdot\|_{B,\zeta  }$ for $\zeta  \in D$ is an hermitian metric defined by the matrix $B(\zeta  )=(B_{i,j}(\zeta  ))_{i,j=1,\ldots, n}$ where $B_{i,j}(\zeta  )=\diffp{^2}{\zeta  _i\partial\overline{\zeta  _j}} \ln b(\zeta  ,\zeta  )$. This means that the Bergman norm of $v=\sum_{i=1}^n v_i e_i$, where $e_1,\ldots, e_n$ is the canonical basis of $\cc^n$, is given by $\|v\|_{B,\zeta  }=\left(\sum_{i,j=1}^n B_{ij}(\zeta  )v_j\overline{v_i}\right)^{\frac12}$.\\
Using Theorems 3.4 and 5.2 of \cite{McN1} and Proposition \ref{prop0.1}, we easily get
\begin{theorem}\mlabel{th3}
For all $\zeta  \in D$ in a neighborhood of $bD$, we have 
$$b(\zeta  ,\zeta  ) \geqs \frac1{\vol \bigl(\cp_{\dzeta}(\zeta  )\bigr)}.$$
Let $w$ be any orthonormal coordinates system centered at $\zeta  $ and let $v_j$
be the unit vector in the $w_j$-direction. Then we have uniformly with respect to $\zeta  $ :
$$\left|\diffp{^{|\alpha|+|\beta |} b }{w^\alpha \partial \overline{w}^\beta} (\zeta  ,\zeta  )\right|\leqs \frac1{\vol\bigl(\cp_{\dzeta}(z)\bigr) \prod_{j=1}^n \tau(\zeta  ,v_j,\dzeta)^{\alpha_j+\beta_j} }.$$
\end{theorem}
Theorem \ref{th3} yields to the following corollary
\begin{corollary}\mlabel{corI.1.1}
Let $\zeta  \in D$ be a point near $bD$, let $w$ be any orthonormal coordinates system centered at $\zeta  $, let $v_j$
be the unit vector in the $w_j$-direction and let $(B_{ij}^w)_{i,j}$ be the Bergman matrix in the $w$-coordinates. Then we have uniformly with respect to $\zeta  $ :
$$\left|\diffp{^{|\alpha|+|\beta|} B^w_{ij}
}{w^\alpha \partial \overline{w}^\beta} (\zeta  )\right|\leqs \frac1{\prod_{j=1}^n \tau(\zeta  ,v_j,\dzeta)^{\alpha_j+\beta_j} }.$$
\end{corollary}
McNeal proved in \cite{McN2} :
\begin{proposition}\mlabel{propI.1.2}
Let $\zeta  \in D$ be a point near $bD$ and let $\lambda_1(\zeta )\geq \lambda_2(\zeta )\geq\ldots\geq \lambda _n(\zeta  )$ be the eigenvalues of $B(\zeta)$. Then uniformly with respect to $\zeta$
$$\lambda _1(\zeta  )\eqs \tau _1(\zeta,\dzeta)^{-2},\
\lambda _2(\zeta  )\eqs \tau _{n}(\zeta,\dzeta)^{-2},
\ldots,\lambda _n(\zeta  )\eqs \tau_2(\zeta  ,\dzeta)^{-2}.$$
\end{proposition}
This also implies that $\det B(\zeta  )\eqs \left(\vol (\cp_{\dzeta}(\zeta  ))\right)^{-1}$, uniformly with respect to $z$.\\
We denote by $e_j(\zeta)$ the $j$-th column of the matrix $A(\zeta)$ so that $e_1(\zeta),\ldots,e_n(\zeta)$ is an orthonormal basis of $\cc^n$ for the Bergman metric. We then end this section with the following proposition :
\begin{proposition}\mlabel{propI.1.11}
 For all vectors fields $u$ and $v$, all smooth $(1,1)$-current $\theta$ and all $\zeta\in D$, we have
 $$\frac{|\theta(\zeta)|[u(\zeta),v(\zeta)]}{k(\zeta,u(\zeta))k(\zeta,v(\zeta))}
 \leqs \sum_{j,k=1}^n \frac{|\theta(\zeta)|[e_j(\zeta),e_k(\zeta)]}{k(\zeta,e_j(\zeta))k(\zeta,e_k(\zeta))}.$$
 \end{proposition}
{\it Proof :} We write $u=\sum_{j=1}^nu_je_j$ and $v=\sum_{j=1}^nv_je_j$. We thus have 
$$\frac{|\theta(\zeta)|[u(\zeta),v(\zeta)]}{k(\zeta,u(\zeta))k(\zeta,v(\zeta))}
 \leqs \sum_{j,k=1}^n |u_j(\zeta)||v_k(\zeta)|\frac{|\theta(\zeta)|[e_j(\zeta),e_k(\zeta)]}{k(\zeta,u(\zeta))k(\zeta,v(\zeta))}.$$
 From Proposition \ref{propI.1.7} we have $k(\zeta,u(\zeta))\eqs\dzeta\|u\|_{B,\zeta}$ and since $e_1(\zeta),\ldots,e_n(\zeta)$ is an orthonormal basis for the Bergman metric $k(\zeta,u(\zeta))\eqs \dzeta\sqrt{\sum_{j=1}^n |u_j(\zeta )|^2}\geqs \dzeta|u_j(\zeta)|$. The same holds true for $v$ and so
$$\frac{|\theta(\zeta)|[u(\zeta),v(\zeta)]}{k(\zeta,u(\zeta))k(\zeta,v(\zeta))}
 \leqs \sum_{j,k=1}^n \frac{|\theta(\zeta)|[e_j(\zeta),e_k(\zeta)]}{\dzeta^2}.$$
Finally, again since $e_1(\zeta),\ldots,e_n(\zeta)$ is an orthonormal basis for the Bergman metric,
$k(\zeta,e_j(\zeta))\eqs\dzeta\|e_j(\zeta)\|_{B,\zeta} \eqs \dzeta$ and the proof of the proposition is complete.\qed
\subsection{Derivatives of the matrix $A$}\mlabel{ssI.2}
We will need upper bounds of $k(\zeta  ,dA_\zeta[u]\cdot\Lambda )$ for any vector $u$. Since
$k(\zeta,dA_\zeta[u]\cdot\Lambda)\eqs \deltar (\zeta)\|dA_\zeta[u]\cdot\Lambda\|_{B,\zeta}$, we look for an upper bound of $\|dA_\zeta[u]\cdot\Lambda\|_{B,\zeta}$.\\
\medskip
First we compute $dA_\zeta$. Let $\ch_n$ be the set of hermitian matrices in $\cc^n$, let $\ch^{++}_n$ be the set of positive definite hermitian matrices in $\cc^n$ and let $\Phi:\ch^{++}_n\to\ch^{++}_n$ be defined by $\Phi(M)=M^{-2}$. The map $\Phi$ is one to one and for all $\zeta  \in D$, $A(\zeta  )=\Phi^{-1}(B(\zeta  ))$. We use the inverse mapping theorem in order to deduce from Corollary \ref{corI.1.1} the needed estimates on $A(\zeta  )$. For $M\in \ch^{++}_n$ and $H\in \ch_n$. We have :
\begin{align*}
d\Phi_M(H)&=M^{-2} (MH+HM)M^{-2}. 
\end{align*}
We want to compute the inverse of $d\Phi_M$. We first notice that $d\Phi_M(H)=H'$ if and only if $MH+HM=M^2H'M^2$. We use ideas of Rosenblum \cite{Ros} in order to solve explicitly this equation. The computations are quiet similar but not exactly the same. We give them for completness.\\
Let $\Omega$ be a bounded open set in $\cc$ such that $\sp(M)$, the spectrum of $M$, is included in $\Omega$, and $\sp(-M)\cap\overline \Omega=\emptyset$. This is always possible because $\sp(M)$ is included in $]0,+\infty[$. We denote by $I_n$ the identity matrix in $\cc^n$. No $\xi$ in $b\Omega$ belongs to $\sp(M)\cup\sp(-M)$, so $\xi I_n+ M$ and $\xi I_n-M$ are invertible and Dunford's functional calculus asserts that
\begin{align}
 \frac1{2i\pi}\int_{b\Omega} (M-\xi I_n)^{-1} M^2H'M^2d\xi&=M^2H'M^2,\mlabel{eqI.1}\\
\frac1{2i\pi}\int_{b\Omega}  M^2H'M^2(M+\xi I_n)^{-1}d\xi&=0.\mlabel{eqI.2}
\end{align}
Therefore, setting $H=\frac1{2i\pi}\int_{b\Omega} (M-\xi I_n)^{-1} M^2H'M^2 (M+\xi I_n)^{-1}d\xi$, Equalities (\ref{eqI.1}) and (\ref{eqI.2}) yield 
\begin{align*}
 MH+HM=&
    \frac1{2i\pi}\int_{b\Omega} M(M-\xi I_n)^{-1} M^2H'M^2 (M+\xi I_n)^{-1}d\xi\\
    &+\frac1{2i\pi}\int_{b\Omega} (M-\xi I_n)^{-1} M^2H'M^2 (M+\xi I_n)^{-1}Md\xi\\
 =&
    \frac1{2i\pi}\int_{b\Omega} (M-\xi I_n)(M-\xi I_n)^{-1} M^2H'M^2 (M+\xi I_n)^{-1}d\xi\\
   & +\frac1{2i\pi}\int_{b\Omega} (M-\xi I_n)^{-1} M^2H'M^2 (M+\xi I_n)^{-1}(M+\xi I_n)d\xi\\
 =& M^2H'M^2.
 \end{align*}
Thus, by the inverse mapping theorem, we get
\begin{align*}
 &d\Phi_{B(\zeta  )}^{-1}(H')\\
 &=
 \frac1{2i\pi}\int_{b\Omega_\zeta} \bigl(\Phi^{-1}({B(\zeta  )})-\xi I_n\bigr)^{-1} \bigl(\Phi^{-1}({B(\zeta  )})\bigr)^2H'\bigl(\Phi^{-1}({B(\zeta  )})\bigr)^2 \bigl(\Phi^{-1}({B(\zeta  )})+\xi I_n\bigr)^{-1}d\xi
\end{align*}
where $\Omega_\zeta$ is any bounded open set in $\cc$ such that $\sp\bigl(\Phi^{-1}({B(\zeta  )})\bigr)$ is included in $\Omega_\zeta$ and $\sp\bigl(-\Phi^{-1}({B(\zeta  )})\bigr)\cap\overline{\Omega_\zeta}=\emptyset.$\\
Let $u$ be a unit vector in $\cc^n$. 
We fix $\zeta_0\in D$ and an orthonormal basis $w_1,\ldots,w_n$ on $\cc^n$, orthogonal for $B(\zeta_0)$.
We denote by $B^w(\zeta)$ the matrix of the Bergman metric in the basis $w_1,\ldots, w_n$ and we assume that
$B^w(\zeta_0)=\begin{pmatrix}
\lambda_1(\zeta_0) &&0\\
&\ddots&\\
0&&\lambda_n(\zeta_0  )                                                                                                                                                                        \end{pmatrix}
$, $\lambda_1(\zeta_0 )> \lambda_2(\zeta_0 )\geq\ldots\geq \lambda _n(\zeta_0  )$.
We denote by $P$ the unitary matrix such that $B(\zeta  _0)=PB^w(\zeta_0)\overline P^t$. We also define the two diagonal matrices
$B^w(\zeta_0)^{-\frac12}=
\begin{pmatrix}
\lambda_1(\zeta_0)^{-\frac12}&&0 \\
&\ddots&\\
0&&\lambda_n(\zeta  )^{-\frac12}                                                                                                                                                                        \end{pmatrix}
$ and 
$B^w(\zeta_0)^{\frac12}=
\begin{pmatrix}
\lambda_1(\zeta_0)^{\frac12} &&0\\
&\ddots&\\
0&&\lambda_n(\zeta  )^{\frac12}                                                                                                                                                                        \end{pmatrix}$ 
so that $A(\zeta  _0)=PB^w(\zeta_0)^{-\frac12}\overline{P}^t.$ \\
We have
$dA_{\zeta_0}[u]\cdot\Lambda=\frac1{2i\pi} P D(\zeta_0)\overline P^t\cdot \Lambda,$ where
\begin{align*}
  D(\zeta_0)=
 \int_{b\Omega_{\zeta_0}} 
 \bigl(B^w(\zeta_0)^{-\frac12}-\xi I_n\bigr)^{-1}
 B^w(\zeta_0)^{-1} \diffp{B^w}{u}(\zeta_0)B^w(\zeta_0)^{-1}
 \bigl(B^w(\zeta_0)^{-\frac12}+\xi I_n\bigr)^{-1} d\xi.
 \end{align*}
Since $P$ is a unitary matrix, we have
\begin{align*}
 \|dA_{\zeta  _0}[u]\cdot\Lambda\|_{B,\zeta_0}&=|B^w(\zeta_0)^{\frac12}\cdot D(\zeta_0)\cdot\overline{P}^t\cdot \Lambda|\\
 &\leq \|B^w(\zeta_0)^{\frac12}D(\zeta_0)\|_\infty\cdot|\Lambda|.
\end{align*}
Corollary \ref{corI.1.1} and Proposition \ref{propI.1.2} imply that $\mu_{i,j}=\diffp{B^w_{i,j}}{u}(\zeta_0,\zeta_0)$ satisfies
\begin{align}
 |\mu_{i,j}|&\leqs \frac{\lambda_i(\zeta  _0)^{\frac12}\lambda _j(\zeta  _0)^{\frac12}}{\tau(\zeta_0,u,\dzetaO )}.\mlabel{eq11}
\end{align}
For $\xi\in b\Omega_{\zeta_0}$, we have
\begin{align}
 &B^w(\zeta_0)^{\frac12}\bigl(B^w(\zeta_0)^{-\frac12}-\xi I_n\bigr)^{-1}
 B^w(\zeta_0)^{-1} \diffp{B^w}{u}(\zeta_0)B^w(\zeta_0)^{-1}
 \bigl(B^w(\zeta_0)^{-\frac12}+\xi I_n\bigr)^{-1}\nonumber\\
 &=\Bigl(\lambda_k(\zeta_0)^{-\frac12} (\lambda_k(\zeta_0)^{-\frac12}-\xi)^{-1}\mu_{kl}\lambda_l(\zeta_0)^{-1}
 (\lambda_l(\zeta_0)^{-\frac12}+\xi)^{-1}\Bigr)_{k,l}.\mlabel{eqI.3'}
 \end{align}
In order to estimate $\|B^w(\zeta_0)^{\frac12}D(\zeta_0)\|_\infty$, we integrate (\ref{eqI.3'}) over $b\Omega_{\zeta_0}$, but before, we choose a good open set $\Omega_{\zeta_0}$. From proposition \ref{propI.1.2}, we have $\lambda_n^{-\frac12}(\zeta_0)\geq\ldots\geq\lambda_1^{-\frac12}(\zeta_0)\eqs\dzetaO $, so there exists $c>0$ sufficiently small such that $\Omega_{\zeta_0}=\cup_{j=1}^n\Delta(\lambda_j^{-\frac12}(\zeta_0),c\dzetaO )$ is included in $\{\xi\in\cc,\ \Re\xi\geq \frac c2\dzetaO \}$. Thus $\sp(\Phi^{-1}(B(\zeta_0)))$ is included in $\Omega_{\zeta_0}$ and $\sp(-\Phi^{-1}({B(\zeta_0  )}))\cap\overline{\Omega_{\zeta_0}}=\emptyset$.\\
For all $k$ and all $\xi\in b\Omega_{\zeta_0}$ holds:
\begin{align}
 |\lambda_k(\zeta_0)^{-\frac12} - \xi |^{-1}&\leqs \dzetaO ^{-1}.\mlabel{eqI.3}
\end{align}
Since $\xi$ belongs to $b\Omega_{\zeta_0}$, there exists $j$ and $\phi$ such that $\xi=\lambda_j(\zeta_0)^{-\frac12}+c\dzetaO e^{i\phi}$, so for all $l$ 
\begin{align*}
 |\lambda_l(\zeta_0)^{-\frac12}+\xi|&\geq \lambda_l(\zeta_0)^{-\frac12} +\lambda_j(\zeta_0)^{-\frac12}-c\dzetaO \nonumber\\
 &\geq \lambda_l(\zeta_0)^{-\frac12}.
\end{align*}
so
\begin{align}
 |\lambda_l(\zeta_0)^{-\frac12}+\xi|^{-1}&\leq \lambda_l(\zeta_0)^{\frac12}.\mlabel{eqI.4}
\end{align}
Inequalities (\ref{eq11}), (\ref{eqI.3}) and (\ref{eqI.4}) yield for all $k$ and $l$ :
\begin{align*}
 &\lambda_k(\zeta_0)^{-\frac12} (\lambda_k(\zeta_0)^{-\frac12}-\xi)^{-1}\mu_{kl}\lambda_l(\zeta_0)^{-1}
 (\lambda_l(\zeta_0)^{-\frac12}+\xi)^{-1}\\
 &\hskip 100pt\leqs 
 \lambda_k(\zeta_0)^{-\frac12}
 \dzetaO ^{-1}
 \frac{\lambda_k(\zeta_0)^{\frac12}\lambda_l(\zeta_0)^{\frac12}}{\tau(\zeta_0,u,\dzetaO )} \lambda_l(\zeta_0)^{-1}\lambda_l(\zeta_0)^{\frac12}\\
 &\hskip 100pt\leqs \frac1{\dzetaO \tau(\zeta_0,u,\dzetaO )}.
\end{align*}
We now integrate the previous inequality on $b\Omega_{\zeta_0}$ and get, since the 
 the length of $b\Omega_{\zeta_0}$ is less than $2\pi n c \dzetaO $ :
 $$\|B^w(\zeta_0)^{\frac12} D(\zeta_0)\|_\infty\leqs \frac1{\tau(\zeta_0,u,\dzetaO )}.$$ 
Thus we have proved the following lemma :
\begin{lemma}\mlabel{lemI.2.1}
 For all $\zeta\in D$ close enough to $bD$, all vector $u\in\cc^n$, all $\Lambda\in\Delta_n(1)$, we have uniformly with respect to $\zeta$, $u$ and $\Lambda$
 $$\|dA_\zeta[u]\cdot\Lambda\|_{B,\zeta}\leqs \frac1{\tau(\zeta,u,\dzeta)}.$$
\end{lemma}
We deduce from Lemma  \ref{lemI.2.1} we the following corollary which will be very useful.
\begin{corollary}\mlabel{lemI.2.2}
 There exists $c>0$ sufficiently small such that for all $\xi\in D$ close to $bD$, all $\zeta\in c\cp_{\dxi}(\xi)$, all $\Lambda\in\Delta_n(1)$, all vector $v\in\cc^n$, we have uniformly
 \begin{align*}
  k(\zeta,A(\xi)\cdot \Lambda)&\leqs \dxi,\\
  k(\zeta,dA_\xi[v]\cdot \Lambda)&\leqs k(\xi,v).\\
 \end{align*}
\end{corollary}
{\it Proof :} From Corollary \ref{corI.1.3} we have $k(\zeta,A(\xi)\cdot \Lambda)\eqs k(\xi,A(\xi)\cdot \Lambda)$ and since $\|A(\xi)\cdot\Lambda\|_{B,\xi}=|\Lambda|$, we get 
\begin{align*}
 k(\zeta,A(\xi)\cdot\Lambda)
 &\eqs \dxi\|A(\xi)\cdot\Lambda\|_{B,\xi}\\
 &\leqs \dxi.
\end{align*}
In the same way we have $k(\zeta,dA_\xi[v]\cdot \Lambda)\leqs \dxi\|dA_\xi[v]\cdot\Lambda\|_{B,\xi}$ and Lemma \ref{lemI.2.1} yields
\begin{align*}
 k(\zeta,dA_\xi[v]\cdot \Lambda)&\leqs \frac{\dxi}{\tau(\xi,v,\dxi)}=k(\xi,v).
\end{align*}
\qed
\subsection{Case $1-t\leq \frac\gamma2\dz$}\mlabel{ssI.3}
In this subsection, we want to estimate 
\begin{align*}
\i&:= \int_{\over{z\in\cp_\varepsilon (z_0)\cap D}{\over{\Lambda\in \Delta_n(\rho)}{t\in[1-\frac\gamma2\dz ,1]}} } \frac{|\theta(\hl)|\left[\diffp{h_\Lambda}t(z,t),d_z\hl[u(z)]\right]}{k(z,u(z))\vol (\Delta_n(\rho))}
dtd\Lambda d\lambda(z).
\end{align*}
We first prove the following lemma  for $1-t\leq \gamma \dz $ and not only for $1-t\leq\frac\gamma2\dz$ because we will also use it in Subsection \ref{ssI.5}.
\begin{lemma}\mlabel{lemI.3.0}
 There exists $C>0$ such that for all $z\in D$ close to $bD$, all $t\in[t_0,1-\gamma\dz]$, all $\Lambda\in\Delta_n(1)$, the point $p=tz+\frac{t(1-t)}{\dz }A(z)\cdot\Lambda$ belongs to $C\frac{1-t}{\dz }\cp_{\dz}(z)$.
\end{lemma}
{\it Proof :} We write $p$ as $p=z-\frac{1-t}{\dz } \dz z+\frac{t(1-t)}{\dz }A(t)\cdot\Lambda$.\\
In the one hand, $z- \dz z$ belongs to $K\cp_\dz(z)$ for some uniform $K$ because $|z-d(z)z-z|\leqs d(z)$.\\
In the other hand, since $\|A(z)\cdot\Lambda\|_{B,z}\leqs 1$,  there exists a uniform $K'$ such that $z+tA(z)\cdot\Lambda$ belongs to $K'\cp_{\dz}(z)$. Therefore, putting $C=K+K'$, $z-\dz(z)+tA(z)\cdot\Lambda$ belongs to $C\cp_{\dz}(z)$ and so $p$ belongs to $C\frac{1-t}{\dz }\cp_{\dz}(z)$.\qed

This lemma gives us the following inequalities :
\begin{corollary}\mlabel{corI.3.2}
 If $\gamma>0$ is small enough, for all $z\in D$ close to $bD$, all $t\in[1-\frac\gamma2\dz,1]$ and all $\Lambda\in\Delta_n(1)$, the following estimates hold :
 \begin{align*}
  \dhl&\eqs \dtz\eqs \dz\\
  \tau(z,v,\dz)&\eqs \tau(tz,v,\dtz)\eqs \tau(h_\Lambda(z,t),v,\dhl),\\
  k(\hl,v)&\eqs k(tz,v)\eqs k(z,v).
 \end{align*}
\end{corollary}
{\it Proof :} Lemma \ref{lemI.3.0} implies that $\hl$ belongs to $c\cp_\dz(z)$, $c$ arbitrary small provided $\gamma$ is small enough. Proposition \ref{propI.1.3} then implies that $\dhl)\eqs\dz$ and Corollary \ref{corI.1.3} implies that $k(\hl,v)\eqs k(z,v)$.\\
Since $|z-tz|\leqs \gamma \dz$, if $\gamma$ is small enough, $tz$ belongs to $c\cp_\dz(z)$ and, in the same way, we have 
$\dz\eqs \dtz$ and $ k(tz,v)\eqs k(z,v)$.\qed
\begin{lemma}\mlabel{lemI.3.3}
Let $c$ be a positive number. If $\gamma>0$ is small enough, for all $z\in D$ close to $bD$, all $t\in[1-\frac\gamma2\dz,1]$ and all $\Lambda\in\Delta_n(1)$, the point $\hl$ belongs to $c\cp_\dz(z)$ and uniformly
 \begin{align*}
 k\left(\hl,\diffp{h_\Lambda}t(z,t) \right)&\leqs 1,\\
 k(\hl,d_zh_\Lambda(z,t)[u])&\leqs k(z,u).
 \end{align*}
\end{lemma}
{\it Proof : } 
From Lemma \ref{lemI.3.0}, since $1-t\leq \frac\gamma2 d(z)$, if $\gamma$ is small enough then $\hl$ belongs to $c\cp_{\dz}(z)$, $c$ arbitrary small, provided that $\gamma$ is small enough.\\
Since $1-t\leq \frac\gamma2\dz$, $\diffp{h_\Lambda}t(z,t)=z+\frac{1-2t}{\dz}A(z)\cdot\Lambda$. 
Propositions \ref{prop0.1} and \ref{propI.1.6} give
\begin{align}
 k(\hl,z)&\leqs |z|\leqs 1\mlabel{eqI.3.1}
\end{align}
Corollary \ref{corI.3.2} implies that $k\left(\hl,\frac{1-2t}{\dz}A(z)\cdot\Lambda\right)\leqs \frac{1-2t}{\dz}k(z,A(z)\cdot\Lambda)$ and  Corollary \ref{lemI.2.2} then gives 
\begin{align*}
k\left(\hl,\frac{1-2t}{\dz}A(z)\cdot\Lambda\right)&\leqs \frac{|1-2t|}{\dz}\dz\leq 1. 
\end{align*}
 With (\ref{eqI.3.1}), this proves the first inequality.\\
For the second one, we have $d_z\hl[u]=tu+t\frac{(1-t)}{\dz^2}\diffp{d}u(z)A(z)\cdot\Lambda +\frac{t(1-t)}\dz dA_z[u]\cdot\Lambda$. Corollary \ref{corI.3.2} gives
\begin{align}
 k(\hl,tu)&\leqs k(z,u).\mlabel{eqI.3.3}
\end{align}
Proposition \ref{propI.1.8} gives $\left|\frac{t(1-t)}{\dz^2}\diffp\deltar u(z)\right|\leqs\frac1{\tau(z,u,\dz)}$ ; Corollary \ref{lemI.2.2} yields $k(\hl,A(z)\cdot\Lambda)\leqs \dz$ and so
\begin{align}
 k\left(\hl,\frac{t(1-t)}{\dz^2}\diffp\deltar u(z) A(z)\cdot\Lambda\right)&\leqs k(z,u).\mlabel{eqI.3.4}
\end{align}
Finally, again Corollary \ref{lemI.2.2} gives $k\left(\hl, \frac{t(1-t)}\dz dA_z[u]\cdot\Lambda\right)\leqs k(z,u)$. With  Inequalities (\ref{eqI.3.3}) and (\ref{eqI.3.4}), it then comes $k(\hl,d_zh_\Lambda(z,t)[u])\leqs k(z,u)$.\qed
\par\bigskip
We now estimate $\i$. From Lemma \ref{lemI.3.3}, it comes 
\begin{align*}
 \i&\leqs
 \int_{\over{z\in\cp_\varepsilon (z_0)\cap D}{\over{\Lambda\in \Delta_n(\rho)}{t\in[1-\frac\gamma2\dz ,z]}} } \frac{|\theta(\hl)|\left[\diffp{h_\Lambda}t(z,t),d_z\hl[u(z)]\right]}{k\bigl(\hl,d_z\hl[u(z)]\bigr)\cdot k\left(\hl,\diffp{h_\Lambda}t(z,t)\right)}dtd\Lambda d\lambda(z).
\end{align*}
Then Proposition \ref{propI.1.11} gives $\i\leqs\sum_{j,k=1}^n \ijk$ where
\begin{align*}
 \ijk&:=
 \int_{\over{z\in\cp_\varepsilon (z_0)\cap D}{\over{\Lambda\in \Delta_n(\rho)}{t\in[1-\frac\gamma2\dz ,1]}} } \frac{|\theta(\hl)|[e_j(\hl),e_k(\hl)]}{k(\hl,e_j(\hl))\cdot k(\hl,e_k(\hl))}dtd\Lambda d\lambda(z).
\end{align*}
For fixed $z$ and $t$, we make the substitution $\zeta=\hl$, $\Lambda$ running over $\Delta_n(\rho)$. From Lemma \ref{lemI.3.0}, when $|\Lambda|\leq \rho,$ the point  $\hl$ belongs to $C\frac{1-t}{\dz} \cp_\dz(z)$ for some big $C>0$. Moreover, $\det_\rr d_\Lambda h(z,t)\eqs \left(\frac{1-t}{\dz}\right)^{2n}(\det_\cc A(z))^2$ and Proposition \ref{propI.1.2} then yields 
$\det_\rr d_\Lambda h(z,t)\eqs \left(\frac{1-t}{\dz}\right)^{2n} \vol(\cp_\dz(z))$. Therefore
\begin{align*}
 \ijk&\leqs
 \int_{\over{z\in\cp_\varepsilon (z_0)\cap D}{\over{t\in[1-\frac\gamma2\dz ,1]}{\zeta\in C\frac{1-t}{\dz}\cp_\dz(z)}} }
 \left(\frac{\dz}{1-t}\right)^{2n} \frac{|\theta(\zeta)|[e_j(\zeta),e_k(\zeta)]}{\vol(\cp_\dz(z))\cdot \ k(\zeta,e_j(\zeta))\cdot k(\zeta,e_k(\zeta))}d\lambda(\zeta)dt d\lambda(z).
\end{align*}
Now we want to change the order of integration. When $z$ belongs to $\cp_\varepsilon (z_0)$, $t$ to $[1-\frac\gamma2\dz,1]$ and $\zeta$ to $C\frac {1-t}{\dz}\cp_{\dz}(z)$, we have $\dz\eqs\dzeta$ provided $\gamma$ is small enough. Thus there exists $K$, not depending on $\gamma$, $z$, $\zeta$ or $t$ such that $1-K\gamma\dzeta\leq t\leq 1$. Therefore if $\gamma$ is small enough, $t$ belongs to $[1-\frac12\dzeta,1]$\\
Because $\delta$ is a pseudodistance, we also have $\delta(\zeta,z)\leqs \delta(\zeta,z)+\delta(z,z_0)\leqs \varepsilon,$ thus, there exists $K'>0$, big enough, such that $\zeta$ belongs to $K'\cp_{\varepsilon}(z_0)$.\\
Since $\zeta$ belongs to $C\frac{1-t}{\dz}\cp_\dz(z)$, we can write $\zeta=z+C\frac{1-t}{\dz}\mu v$ with $\mu\in\cc$, $v\in\cc^n$, $|v|=1$, such that $|\mu|<\tau(z,v,\dz)$. Provided $\gamma$ is small enough, we have $\dz\eqs\dzeta$ and $\tau(z,v,\dz)\eqs\tau(\zeta,v,\dzeta)$. Therefore $z=\zeta-\frac{1-t}{\dz}\mu v$ with $|\mu|\leqs \tau(\zeta,v,\dzeta)$ and there exists $K''>0$ big enough, such that $z$ belongs to $K''\frac{1-t}{\dzeta}\cp_{\dzeta}(\zeta  )$.\\
Therefore, the set 
$\{(z,t,\zeta),\ z\in\cp_\varepsilon (z_0),\ t\in[1-\frac\gamma2\dz ,1],\ \zeta\in C\frac{1-t}{\dz } \cp_\dz(z)\}$ is included in 
$\{(z,t,\zeta),\ \zeta\in K'\cp_{\varepsilon} (z_0),\ t\in[1-\frac12\dzeta,1],\ z\in K''\frac{1-t}{\dzeta} \cp_{\dzeta}(\zeta)\}$.\\
Moreover, $\vol(\cp_\dz(z))\eqs\vol( \cp_{\dzeta}(\zeta))$ which gives
\begin{align*}
 \ijk&\leqs
 \int_{\over{\zeta\in K'\cp_\varepsilon (z_0)\cap D}{\over{t\in[1-\frac12\dzeta,1]}{z\in K''\frac{1-t}{\dzeta}\cp_{\dzeta}(\zeta)}} }
 \left(\frac{\dzeta}{1-t}\right)^{2n} \frac{|\theta(\zeta)|[e_j(\zeta),e_k(\zeta)]}{\vol(\cp_{\dzeta}(\zeta))\cdot \ k(\zeta,e_j(\zeta))\cdot k(\zeta,e_k(\zeta))}d\lambda(z)dt d\lambda(\zeta).
\end{align*}
Integrating for $z$ in $K''\frac{1-t}{\dzeta}\cp_{\dzeta}(\zeta)$ we get
\begin{align*}
 \ijk&\leqs
 \int_{\over{\zeta\in K'\cp_\varepsilon (z_0)\cap D}{t\in[1-\frac12\dzeta,1]} }
  \frac{|\theta(\zeta)|[e_j(\zeta),e_k(\zeta)]}{ k(\zeta,e_j(\zeta))\cdot k(\zeta,e_k(\zeta))}dt d\lambda(\zeta).
\end{align*}
Now we integrate for $t\in[1-\frac12 \dzeta,1]$ and we obtain
\begin{align*}
 \ijk&\leqs
 \int_{\zeta\in K'\cp_\varepsilon (z_0)\cap D}
  \frac{\dzeta\ |\theta(\zeta)|[e_j(\zeta),e_k(\zeta)]}{ k(\zeta,e_j(\zeta))\cdot k(\zeta,e_k(\zeta))} d\lambda(\zeta)
\end{align*}
and since $\deltap\theta$ is a $(1,1)$-Carleson current, we get
\begin{align*}
 \ijk&\leqs \sigma(\cp_\varepsilon (z_0)\cap bD) \|\deltap\theta\|_{W^1_{1,1}}.
\end{align*}
This finally shows that $\i\leqs\|\deltap\theta\|_{W^1_{1,1}}\sigma(\cp_\varepsilon (z_0)\cap bD).$
\subsection{Case $\gamma\dz\leq 1-t$}\mlabel{ssI.4}
This subsection is devoted to the estimate of
\begin{align*}
\ii&:= \int_{\over{z\in\cp_\varepsilon (z_0)\cap D}{\over{\Lambda\in \Delta_n(\rho)}{t\in[t_0,1-\gamma\dz ]}} } \frac{|\theta(\hl)|\left[\diffp{h_\Lambda}t(z,t),d_z\hl[u(z)]\right]}{k(z,u(z))\vol (\Delta_n(\rho))}
dtd\Lambda d\lambda(z).
\end{align*}
We first look for estimates of $k(\hl,d_zh_\Lambda(z,t)[u])$ and $k(\hl,\diffp{h_\Lambda}t(z,t))$. We begin with following lemma that we prove for $z$ and $t$ such that $\frac\gamma2\dz\leq 1-t$.
\begin{lemma}\mlabel{lemI.4.1}
 There exists $C>0$ such that for all $z\in D$ close to $bD$, all $t\in[t_0,1-\frac\gamma2\dz]$, all $\Lambda\in\Delta_n(\rho)$, the point $q=tz+tA(tz)\cdot\Lambda$ belongs to $Ct\rho \cp_{\dtz}(tz)$.
\end{lemma}
{\it Proof :} If we write $A(tz)\cdot\Lambda$ as $\mu v$ with $\mu\in\cc$ and $v\in\cc^n$, $|v|=1$. We have
$$\frac{|\mu|}{\tau(tz,v,\dtz)}\eqs \|A(tz)\cdot\Lambda\|_{B,tz}\leqs\rho.$$
Thus there exists a uniform $K>0$ such that $|\mu|\leq K\rho\tau(tz,v,\dtz)$ and so, $q$ belongs to $t\rho C\cp_{\dtz}(tz)$ for some $C$ which does not depend on $z$, $t$ or $\rho$. \qed
\begin{corollary}\mlabel{corI.4.3}
 If $\rho>0$ is small enough, for all $z\in D$ close to $bD$, all $t\in[t_0,1-\gamma\dz]$ and all $\Lambda\in\Delta_n(\rho)$ :
 \begin{align*}
  \dhl&\eqs \dtz\eqs 1-t\\
  \tau(z,v,1-t)&\eqs \tau(tz,v,\dtz)\eqs \tau(h_\Lambda(z,t),v,\dhl),\\
  k(\hl,v)&\eqs k(tz,v)\eqs \frac{1-t}{\tau(z,v,1-t)}.
 \end{align*}
\end{corollary}
{\it Proof :}
Firstly, we have $\dtz=1-t+t\dz \eqs1-t$. Secondly, since from Lemma \ref{lemI.4.1}, $h_\Lambda(z,t)$ belongs to $Ct\rho\cp_{\dtz}(tz)$, choosing $\rho$ sufficiently small we get from Proposition \ref{propI.1.3}
$\dhl\eqs\dtz$. This prove the first chain of almost equalities.
\par\medskip
Since $|z-tz|\leqs 1-t$, $z$ belongs to $\cp_{K(1-t)}(tz)$ and since $1-t\eqs \dtz$, from Propositions \ref{propI.1.4} and \ref{propI.1.5}, we have 
$\tau(z,v,1-t)\eqs\tau(tz,v,\dtz)$ and $\frac{1-t}{\tau(z,v,1-t)} \eqs k(tz,v)$.
Now since $h_\Lambda(z,t)$ belongs to $Ct\rho\cp_{\dtz}(tz)$, Proposition  \ref{propI.1.4} gives $\tau(tz,v,\dtz)\eqs\tau(h_\Lambda(z,t),v,\dtz)$, provided $\rho$ is small enough.  Since  $\dhl\eqs \dtz$, it then comes from Proposition \ref{propI.1.5} $\tau(h_\Lambda(z,t),v,\dtz)\eqs\tau(h_\Lambda(z,t),v,\dhl)$ and thus $k(\hl,v)\eqs k(tz,v)$. \qed
\begin{lemma}\mlabel{lemI.4.2}
If $\rho>0$ is small enough, for all $z\in D$ close to $bD$, all $t\in[t_0,1-\gamma\dz]$ and all $\Lambda\in\Delta_n(\rho)$, the following inequalities hold :
 \begin{align*}
 k(\hl,d_zh_\Lambda(z,t)[u])&\leqs k(z,u)\left(\frac{1-t}{\dz }\right)^{1-\frac1m},\\
 k\left(\hl,\diffp{h_\Lambda}t(z,t) \right)&\leqs 1.
 \end{align*}
\end{lemma}
{\it Proof : }
Since $1-t\geq \gamma\dz $, $h_\Lambda(z,t)=tz+tA(tz)\cdot\Lambda$ and thus $d_zh_\Lambda(z,t)[u]=tu+tdA_{tz}[u]\cdot\Lambda$. Lemma \ref{lemI.4.1} implies that $\hl$ belongs to $c\cp_{\dz}(tz)$, $c$ as small as needed if $\rho$ is small enough. We then get from Corollary \ref{corI.4.3}  
\begin{align}
 k(\hl,tu) &\eqs\frac{1-t}{\tau(z,u,1-t)}.\mlabel{eqI.4.1}
\end{align}
Using successively  Corollary \ref{lemI.2.2} and \ref{corI.4.3} we get $k(\hl,tdA_{tz}[u]\cdot\Lambda)\leqs\frac{1-t}{\tau(z,u,1-t)}$.
With (\ref{eqI.4.1}), this yields 
$$k(\hl,d_z\hl[u])\leqs\frac{1-t}{\tau(z,u,1-t)}.$$
Proposition \ref{propI.1.5} then implies
$$k(\hl,d_z\hl[u])\leqs \left(\frac{\dz }{1-t}\right)^{\frac1m} \frac{1-t}{\tau(z,u,\dz )},$$ which proves the first inequality.
\par\medskip
We now prove the second inequality. We have $\diffp{h_\Lambda}t(z,t) =z+A(tz)\cdot\Lambda +tdA_{tz}[z]\cdot\Lambda.$ 
Corollary \ref{corI.4.3} gives
\begin{align}
 k(\hl,z)&\leqs k(tz,z)\mlabel{eqI.4.2}
\end{align}
Next from Corollary \ref{lemI.2.2} we get
\begin{align}
 k(\hl,A(tz)\cdot\Lambda) &\leqs \dz \mlabel{eqI.4.3}
\end{align}
and again with Corollary \ref{lemI.2.2} we have
\begin{align}
 k(\hl,tdA_{tz}[z]\cdot\Lambda) &\leqs k(tz,z).\mlabel{eqI.4.4}
\end{align}
Putting together the inequalities (\ref{eqI.4.2}), (\ref{eqI.4.3}) and (\ref{eqI.4.4}) we obtain
$$k\left(\hl,\diffp{h_\Lambda}t(z,t)\right)\leqs k(tz,z)$$ 
and Propositions \ref{prop0.1} and \ref{propI.1.6} end the proof of the lemma.\qed
\par\medskip
We now come to the heart of the matter of this subsection : We estimate $\ii$. Lemma \ref{lemI.4.2} immediately gives
\begin{align*}
 \ii&\leqs
 \int_{\over{z\in\cp_\varepsilon (z_0)\cap D}{\over{\Lambda\in \Delta_n(\rho)}{t\in[t_0,1-\gamma\dz ]}} }\left(\frac{1-t}{\dz }\right)^{1-\frac1m} \frac{|\theta(\hl)|\left[\diffp{h_\Lambda}t(z,t),d_z\hl[u(z)]\right]dtd\Lambda d\lambda(z)}{k\bigl(\hl,d_z\hl[u(z)]\bigr)\cdot k\left(\hl,\diffp{h_\Lambda}t(z,t)\right)}.
\end{align*}
Then Proposition \ref{propI.1.11} gives $\ii\leqs\sum_{j,k=1}^n \iijk$ where
\begin{align*}
 \iijk&:=
 \int_{\over{z\in\cp_\varepsilon (z_0)\cap D}{\over{\Lambda\in \Delta_n(\rho)}{t\in[t_0,1-\gamma\dz ]}} }\left(\frac{1-t}{\dz }\right)^{1-\frac1m} \frac{|\theta(\hl)|[e_j(\hl),e_k(\hl)]dtd\Lambda d\lambda(z)}{k(\hl,e_j(\hl))\cdot k(\hl,e_k(\hl))}.
\end{align*}
We make the substitution $\zeta=\hl$ for $\Lambda$ running over $\Delta_n(\rho)$. Since, Proposition \ref{propI.1.2}, $\det_\rr A(tz)\eqs\vol(\cp_{\dtz}(tz))$, and, since Lemma \ref{lemI.4.1}, $\{\hl,\ |\Lambda|<\rho\}\subset Ct\rho\cp_\dtz(tz)$, we have
\begin{align*}
 \iijk&\leqs
 \int_{\over{z\in\cp_\varepsilon (z_0)\cap D}{\over{\zeta\in C\rho\cp_\dtz(tz)}{t\in[t_0,1-\gamma\dz ]}} }\left(\frac{1-t}{\dz }\right)^{1-\frac1m}\frac1{\vol(\cp_\dtz(tz))} \frac{|\theta(\zeta  )|[e_j(\zeta  ),e_k(\zeta  )]}{k(\zeta  ,e_j(\zeta  ))k(\zeta  ,e_k(\zeta  ))}d\lambda(\zeta)dt d\lambda(z).
\end{align*}
We split $\iijk$ in two parts : $\iijkd$ where we integrate only for $t\in[1-\varepsilon ,1-\gamma\dz ]$ and $\iijke$ where $t$ runs over $[t_0,1-\varepsilon]$. We begin with the easiest part : $\iijke$.\\
If $\rho$ is small enough, $\dzeta\eqs\dtz\eqs 1-t$, thus
\begin{align*}
 \iijke&\leqs
 \int_{\over{z\in\cp_\varepsilon (z_0)\cap D}{\over{\zeta\in C\rho\cp_\dtz(tz)}{t\in[t_0,1-\varepsilon ]}} }
 \frac{(1-t)^{-\frac1m}\dz ^{\frac1m-1}} {\vol(\cp_\dtz(tz))}
 \frac{\dzeta|\theta(\zeta  )|[e_j(\zeta  ),e_k(\zeta  )]}{k(\zeta  ,e_j(\zeta  ))k(\zeta  ,e_k(\zeta  ))}d\lambda(\zeta) dt d\lambda(z).
\end{align*}
We will estimate the integral in $\zeta$ using the fact that $\deltap\theta$ is a $(1,1)$-Carleson current. For fixed $z$ and $t$, let $\zeta_0\in bD$ be a point such that $tz=\zeta_0+\alpha \eta_{\zeta_0}$, $\alpha\in\rr$. Then $|\alpha|\eqs \dtz$ and from Proposition \ref{propI.1.9}, there exists $K>0$, not depending on $z$ or on $t$ such that $\cp_\dtz(tz)\cap D$ is included in $\cp_{K\dtz}(\zeta_0)\cap D$. Moreover, we have $\vol (\cp_\dtz(tz))\eqs\vol (\cp_{K\dtz}(\zeta_0))$. Thus
\begin{align*}
  &\frac1{\vol (\cp_\dtz(tz))}\int_{\zeta\in C\rho\cp_\dtz(tz)}
 \frac{\dzeta|\theta(\zeta  )|[e_j(\zeta  ),e_k(\zeta  )]}{k(\zeta  ,e_j(\zeta  ))k(\zeta  ,e_k(\zeta  )}d\lambda(\zeta)\\
& \leqs
 \frac1{\vol( \cp_\dtz(\zeta_0))}\int_{\zeta\in K\cp_\dtz(\zeta  _0)\cap D}
 \frac{\dzeta|\theta(\zeta  )|[e_j(\zeta  ),e_k(\zeta  )]}{k(\zeta  ,e_j(\zeta  ))k(\zeta  ,e_k(\zeta  )}d\lambda(\zeta)\\
& \leqs \|\theta\|_{W^1_{1,1}} \dtz^{-1}.
\end{align*}
So with Corollary \ref{corI.4.3} we get
\begin{align*}
 \iijke&\leqs
 \|\theta\|_{W^1_{1,1}} \int_{\over{z\in\cp_\varepsilon (z_0)\cap D}{t\in[t_0,1-\varepsilon ]}} 
 (1-t)^{-\frac1m-1} \dz ^{\frac1m-1}
 dt d\lambda(z)\\
 &\leqs \varepsilon^{-\frac1m}\|\theta\|_{W^1_{1,1}} 
 \int_{z\in\cp_\varepsilon (z_0)\cap D}
 \dz ^{\frac1m-1}
  d\lambda(z)\\
 &\leqs\|\theta\|_{W^1_{1,1}} \sigma(\cp_\varepsilon (z_0)\cap bD).
\end{align*}
\par\medskip
Now we deal with $\iijkd$. We write $z\in\cp_\varepsilon (z_0)$ as $z=sz'$ where $s$ belongs to $[1-\varepsilon,1]$ and $z'$ to $\cp_\varepsilon (z_0)\cap bD$. We then have $\dszp=1-s$ and
\begin{align*}
 \iijkd&\leqs
 \int_{\over{z'\in\cp_\varepsilon (z_0)\cap bD}{\over {s\in[1-\varepsilon ,1]}{\over{t\in[1-\varepsilon,1 ]}{\zeta\in C\rho\cp_\dtsz(tsz')}}}}
 \left(\frac{1-t}{1-s}\right)^{1-\frac1m}  \frac{|\theta(\zeta  )|[e_j(\zeta  ),e_k(\zeta  )]\  d\lambda(\zeta)dtds d\sigma(z')}{\vol(\cp_\dtsz(tsz')) \cdot k(\zeta  ,e_j(\zeta  ))\cdot k(\zeta  ,e_k(\zeta  ))}.
\end{align*}
We make the substitution $r=st$ for $t\in[1-\varepsilon ,1]$ and, since $s$ is far from 0, we get
\begin{align*}
 \iijkd&\leqs
 \int_{\over{z'\in\cp_\varepsilon (z_0)\cap bD}{\over {s\in[1-\varepsilon ,1]}{\over{r\in[s(1-\varepsilon),s ]}{\zeta\in C\rho\cp_\drz(rz')}}}}
 \left(\frac{s-r}{1-s}\right)^{1-\frac1m}  \frac{|\theta(\zeta  )|[e_j(\zeta  ),e_k(\zeta  )]\  d\lambda(\zeta)drds d\sigma(z')}{\vol(\cp_\drz(rz')) \cdot k(\zeta  ,e_j(\zeta  ))\cdot k(\zeta  ,e_k(\zeta ) )}.
 \end{align*}
Changing the order of integration between $r$ and $s$ then yields
 \begin{align*}
 \iijkd 
 &\leqs 
 \int_{\over{z'\in\cp_\varepsilon (z_0)\cap bD}{\over {r\in[(1-\varepsilon)^2 ,1]}{\over{\zeta\in C\rho\cp_\drz(rz')}{s\in[r,1 ]}}}}
 \left(\frac{s-r}{1-s}\right)^{1-\frac1m}  \frac{|\theta(\zeta  )|[e_j(\zeta  ),e_k(\zeta  )]\  dsd\lambda(\zeta) dr d\sigma(z')}{\vol(\cp_\drz(rz')) \cdot k(\zeta  ,e_j(\zeta  ))\cdot k(\zeta  ,e_k(\zeta))}
 \end{align*}
Now, integrating separatly for $s\in \left[r,\frac{r+1}2\right]$ and for $s\in\left[\frac{r+1}2,1\right]$, we easily get 
$$\int_{s\in [r,1]} \left(\frac{s-r}{1-s}\right)^{1-\frac1m}ds\leqs1-r.$$
Therefore
\begin{align*}
 \iijkd &\leqs 
 \int_{\over{z'\in\cp_\varepsilon (z_0)\cap bD}{\over {r\in[(1-\varepsilon)^2 ,1]}{\zeta\in C\rho\cp_\drz(rz')}}}
 (1-r) \frac{|\theta(\zeta  )|[e_j(\zeta  ),e_k(\zeta  )]\  d\lambda(\zeta)dr d\sigma(z')}{\vol(\cp_\drz(rz')) \cdot k(\zeta  ,e_j(\zeta  ))\cdot k(\zeta  ,e_k(\zeta  ))}.
 \end{align*}
 If $\rho$ is sufficiently small, since $\zeta$ belongs to $C\rho \cp_\drz(rz')$, we have $\dzeta\eqs \drz\eqs 1-r$. On the other hand, when $r$ belongs to $[(1-\varepsilon )^2,1]$ and $z'$ to $\cp_\varepsilon (z_0)\cap bD$, $rz'$ belongs to $\cp_{2\varepsilon }(z_0)$ because $1-r\leq 2\varepsilon $. Therefore, putting $z=rz'$, we obtain
\begin{align*}
 \iijkd &\leqs 
 \int_{\over{z\in\cp_{2\varepsilon} (z_0)\cap D}{\zeta\in C\rho\cp_\dz(z)}}
  \frac{\dzeta|\theta(\zeta  )|[e_j(\zeta  ),e_k(\zeta  )]}{\vol(\cp_\dz(z)) \cdot k(\zeta  ,e_j(\zeta  ))\cdot k(\zeta  ,e_k(\zeta  ))} d\lambda(\zeta) d\lambda(z).
 \end{align*}
 Now we want to apply Fubbini's Theorem. We notice that when $\zeta$ belongs to $C\rho\cp_\dz(z)\cap D$ and $z$ to $\cp_{2\varepsilon} (z_0)\cap D$, we have $\delta(\zeta  ,z_0)\leqs \varepsilon $, thus $\zeta  $ belongs $K\cp_\varepsilon (z_0)\cap D$ for some big $K$, uniform with respect to $\zeta,z, z_0$ and $\varepsilon$. Moreover, still because $\zeta$ belongs to $C\rho\cp_\dz(z)$,  if $\rho$ is small enough, $z$ belongs to $\cp_\dz(\zeta)$. We also have, still if $\rho$ is small enough, $\dzeta\eqs\dz$, thus  $\frac1K\cp_{\dzeta}(\zeta)\subset \cp_{\dz}(z)\subset K\cp_{\dzeta}(\zeta)$ for some perhaps bigger $K$. In particular $\vol(\cp_{\dzeta}(\zeta))\eqs\vol (\cp_\dz(z))$ and so 
\begin{align*}
 \iijkd &\leqs 
 \int_{\over{\zeta\in\cp_{K\varepsilon} (z_0)\cap D}{z\in K\cp_{\dzeta}(\zeta)}}
  \frac{\dzeta|\theta(\zeta  )|[e_j(\zeta  ),e_k(\zeta  )]}{\vol(\cp_{\dzeta}(\zeta)) \cdot k(\zeta  ,e_j(\zeta  ))\cdot k(\zeta  ,e_k(\zeta  ))}d\lambda(z)\ d\lambda(\zeta)\\
  &\leqs
  \int_{\zeta\in\cp_{K\varepsilon} (z_0)\cap D}
  \frac{\dzeta|\theta(\zeta  )|[e_j(\zeta  ),e_k(\zeta  )]}{k(\zeta  ,e_j(\zeta  ))\cdot k(\zeta  ,e_k(\zeta  ))}\ d\lambda(\zeta)\\
  &\leqs \sigma(\cp_\varepsilon (z_0)\cap bD)\|\deltap\theta\|_{W^1_{1,1}}.
 \end{align*}
 This finally shows that $\ii\leqs\|\deltap\theta\|_{W^1_{1,1}}\sigma(\cp_\varepsilon (z_0)\cap bD).$
 
 \subsection{Case $\frac \gamma2\dz\leq 1-t\leq \gamma\dz $}\mlabel{ssI.5}
The last piece of $H\theta$ that is left to be estimated is
\begin{align*}
\iii&:= \int_{\over{z\in\cp_\varepsilon (z_0)\cap D}{\over{t\in[1-\gamma\dz ,1-\frac\gamma2\dz]} {\Lambda\in \Delta_n(\rho)}}} \frac{|\theta(\hl)|\left[\diffp{h_\Lambda}t(z,t),d_z\hl[u(z)]\right]}{k(z,u(z))\vol (\Delta_n(\rho))}
d\Lambda dtd\lambda(z).
\end{align*}
As in the previous subsections, we first want upper bounds for $k(\hl,d_zh_\Lambda(z,t)[u])$ and $k(\hl,\diffp{h_\Lambda}t(z,t))$.
\begin{lemma}\mlabel{lemI.5.1}
 Let $c$ be a positive number. If $\gamma>0$ and $\rho>0$ are small enough, for all $z\in D$ close to $bD$, all $t\in[1-\gamma\dz, 1-\frac\gamma2\dz]$ and all $\Lambda\in\Delta_n(\rho)$, the point $\hl$ belongs to $c\cp_\dz(z)$ and 
 \begin{align*}
  \dhl&\eqs \dtz\eqs \dz\\
  \tau(z,v,\dz)&\eqs \tau(tz,v,\dtz)\eqs \tau(h_\Lambda(z,t),v,\dhl),\\
  k(\hl,v)&\eqs k(tz,v)\eqs k(z,v).
 \end{align*}
\end{lemma}
{\it Proof :} Lemma \ref{lemI.3.0} implies that $tz+ \frac{t(1-t)}{\dz} A(z)$ belongs $Ct\frac{1-t}{\dz}\cp_\dz(z)$. Thus, for any arbitrary $c$, if $\gamma$ is small enough, $tz+ \frac{t(1-t)}{\dz} A(z)$ belongs $c\cp_\dz(z)$.\\
Lemma \ref{lemI.4.1} implies that $tz+ tA(tz)\cdot\Lambda$ belongs to $ Ct\rho\cp_\dz(z)$. 
Thus, for any arbitrary $c$, if $\rho$ is small enough, $tz+ tA(tz)\cdot\Lambda$ belongs to  $c\cp_\dz(z)$.\\
By convexity $\hl=tz+ \varphi\left(\frac{1-t}{\gamma\dz}\right) \frac{t(1-t)}{\dz} A(z)+ \left(1-\varphi\left(\frac{1-t}{\gamma\dz}\right)\right)tA(tz)\cdot\Lambda$ belongs to $c\cp_\dz(z)$.\\
The rest of the proof is exactly as in Corollary \ref{corI.3.2}, so we omit it.\qed

\begin{lemma}\mlabel{lemI.5.2}
If $\rho>0$ and $\gamma>0$ are small enough, for all $z\in D$ close to $bD$, all $t\in[1-\gamma\dz,1-\frac\gamma2\dz]$ and all $\Lambda\in\Delta_n(\rho)$, the following inequalities hold :
 \begin{align*}
 k\left(\hl,\diffp{h_\Lambda}t(z,t)\right )&\leqs 1,\\
 k\left(\hl,d_zh_\Lambda(z,t)[u]\right)&\leqs k(z,u).
 \end{align*}
\end{lemma}
{\it Proof : }
We recall that $h_\Lambda(z,t)=tz+ t\varphi \left(\frac{1-t}{\gamma\dz}\right)\frac{1-t}{\dz} A(z) \cdot\Lambda   +t\left(1-\varphi \left(\frac{1-t}{\gamma\dz}\right)\right)A(tz) \cdot\Lambda$ so  
\begin{align*}
 &\diffp{h_\Lambda}t(z,t)=\\
 &= z+ \left(\frac{1-2t}{\dz}\varphi\left(\frac{1-t}{\gamma\dz}\right)-\frac{t(1-t)}{\gamma\dz^2} \varphi'\left(\frac{1-t}{\gamma\dz}\right)\right) A(z)\cdot \Lambda\\
 &+\left(1-\varphi\left(\frac{1-t}{\gamma\dz}\right) +\frac t{\gamma\dz} \varphi'\left(\frac{1-t}{\gamma\dz}\right)\right) A(tz)\cdot\Lambda+ t\left(1-\varphi\left(\frac{1-t}{\gamma\dz}\right)\right) dA_{tz}[z]\cdot\Lambda.
\end{align*}
On the one hand
\begin{align}
 k(\hl,z)&\leqs|z|\leqs 1\mlabel{eqI.5.4}.
\end{align}
On the other hand
\begin{align*}
 &k\left(\hl,\left(\frac{1-2t}{\dz}\varphi\left(\frac{1-t}{\gamma\dz}\right)-\frac{t(1-t)}{\gamma\dz^2} \varphi'\left(\frac{1-t}{\gamma\dz}\right)\right) A(z)\cdot \Lambda\right)\\
 &\leqs \frac1{\dz } k(\hl,A(z)\cdot\Lambda),
\end{align*}
and Corollary \ref{lemI.2.2} and Lemma \ref{lemI.5.1} then give
\begin{align}
 k\left(\hl,\left(\frac{1-2t}{\dz}\varphi\left(\frac{1-t}{\gamma\dz}\right)-\frac{t(1-t)}{\gamma\dz^2} \varphi\left(\frac{1-t}{\gamma\dz}\right)\right) A(z)\cdot \Lambda\right)&\leqs 1.\mlabel{eqI.5.3}
\end{align}
Similarly, Corollary \ref{lemI.2.2} and Lemma \ref{lemI.5.1} give
\begin{align}
k\left(\hl, \left(\hskip-1pt 1-\varphi\left(\frac{1-t}{\gamma\dz}\right) +\frac t{\gamma\dz} \varphi'\left(\frac{1-t}{\gamma\dz}\right)\right) A(tz)\cdot\Lambda\hskip-1pt\right)&\leqs \frac1\dz k(\hl,A(tz)\cdot\Lambda)\nonumber\\
&\leqs 1.\mlabel{eqI.5.2}
\end{align}
Again with Corollary \ref{lemI.2.2} and Lemma \ref{lemI.5.1}, we obtain
\begin{align}
k\left(\hl, t\left(1-\varphi\left(\frac{1-t}{\gamma\dz}\right)\right) dA_{tz}[z]\cdot\Lambda\right)&\leqs  k(\hl,dA_{tz}[z]\cdot\Lambda)\nonumber\\
&\leqs |z|\leqs 1.\mlabel{eqI.5.1}
\end{align}
Together (\ref{eqI.5.4}), (\ref{eqI.5.3}), (\ref{eqI.5.2}) and (\ref{eqI.5.1}) give $k\left(\hl,\diffp{h_\Lambda}t(z,t) \right)\leqs 1$.\\
We now prove the second inequality :
\begin{align*}
 d_z\hl[u]=&tu + 
 \left(-\frac{t(1-t)}{\dz^2}\diffp \deltar u(z)\varphi \left(\frac{1-t}{\gamma \dz}\right)
 -\frac{t(1-t)^2}{\gamma \dz^2} \diffp \deltar u (z) \varphi' \left(\frac{1-t}{\gamma \dz}\right)
 \right) A(z)\cdot\Lambda\\
 & + \frac{t(1-t)}{\dz} \varphi \left(\frac{1-t}{\gamma \dz}\right) dA_z[u]\cdot\Lambda +t\diffp\deltar u(z) \frac{1-t}{\gamma\dz^2} \varphi' \left(\frac{1-t}{\gamma \dz}\right) A(tz)\cdot\Lambda\\
 &+ t\left(1-\varphi \left(\frac{1-t}{\gamma \dz}\right)\right)dA_{tz}[u]\cdot\Lambda.
\end{align*}
Propisition \ref{propI.1.8} implies that $\left|\diffp\deltar u(z)\right|\leqs\frac{\dz }{\tau(z,u,\dz )}$. Therefore, since $1-t\eqs \dz$, we have
\begin{align*}
 \left|\frac{t(1-t)}{\dz^2}\diffp \deltar u(z)\varphi \left(\frac{1-t}{\gamma \dz}\right)
 -\frac{t(1-t)^2}{\gamma \dz^2} \diffp \deltar u (z) \varphi' \left(\frac{1-t}{\gamma \dz}\right)
 \right|&\leqs \frac1{\tau(z,u,\dz)};\\
 \left|t\diffp\deltar u(z) \frac{1-t}{\gamma\dz^2} \varphi' \left(\frac{1-t}{\gamma \dz}\right)\right|&\leqs  \frac1{\tau(z,u,\dz)}
 \end{align*}
We then get with Corollary \ref{lemI.2.2} and Lemma \ref{lemI.5.1}
\begin{align*}
k(\hl, d_z\hl[u])&\leqs k(z,u).
\end{align*}
\qed
\par\medskip
We now estimate $\iii$. The way is essentially the same as in the previous subsections, the main difference being when we substitute $\zeta=\hl$. By Lemma \ref{lemI.5.2} we get
\begin{align*}
 \iii&\leqs
 \int_{\over{z\in\cp_\varepsilon (z_0)\cap D}{\over{t\in[1-\gamma\dz ,1-\frac\gamma2\dz]}{\Lambda\in \Delta_n(\rho)}}} \frac{|\theta(\hl)|\left[\diffp{h_\Lambda}t(z,t),d_z\hl[u(z)]\right]}{k\bigl(\hl,d_z\hl[u(z)]\bigr)\cdot k\left(\hl,\diffp{h_\Lambda}t(z,t)\right)}d\Lambda dt d\lambda(z).
\end{align*}
As previously, Proposition \ref{propI.1.11} gives $\iii\leqs\sum_{j,k=1}^n \iiijk$ where
\begin{align*}
 \iiijk&:=
 \int_{\over{z\in\cp_\varepsilon (z_0)\cap D}{{\over{t\in[1-\gamma\dz ,1-\frac\gamma2\dz]}{\Lambda\in \Delta_n(\rho)} }}} \frac{|\theta(\hl)|[e_j(\hl),e_k(\hl)]}{k(\hl,e_j(\hl))\cdot k(\hl,e_k(\hl))}d\Lambda dt d\lambda(z).
\end{align*}
Now we make the substitution $\zeta=\hl$ for $\Lambda$ running over $\Delta_n(\rho)$. By Lemma \ref{lemI.5.1} $\hl$ belongs to $\cp_\dz(z)$ and $\dhl)\eqs \dz$. We have to be a little careful with the determinant of the Jacobian matrix of $\hl$. We have
\begin{align*}
\det_\rr (d_\Lambda \hl) &=\left|\det_\cc \left(\varphi\left(\frac{1-t}{\gamma\dz}\right)\frac{1-t}{\dz} A(t) +\left(1-\varphi\left(\frac{1-t}{\gamma\dz} \right)\right)A(tz) \right)\right|^2.
\end{align*}
Since $\frac{1-t}\dz A(t)$ and $A(tz)$ are both positive definite hermitian matrices, we have
\begin{align*}
\det_\rr (d_\Lambda \hl) &\geq
\left|\det_\cc \left(\frac{1-t}{\dz} A(t)\right)\right|^{2\phi}
\left|\det_\cc\left(A(tz) \right)\right|^{2(1-\phi)},
\end{align*}
where $\phi$ is a shortcut for $\varphi\left(\frac{1-t}{\gamma\dz}\right)$.\\
Since $\frac{1-t}\dz \eqs1$, Proposition \ref{propI.1.2} gives  $\det_\cc \left(\frac{1-t}{\dz} A(t)\right)\eqs\bigl(\vol (\cp_\dz(z))\bigr)^{\frac12}$, uniformly with respect to $z$.\\
Again using Proposition \ref{propI.1.2}, we get $\det_\cc A(tz)\eqs \bigl(\vol (\cp_\dtz(tz))\bigr)^{\frac12}$. Since $tz$ belongs to $\cp_{K\dz}(z)$ for some uniform big $K$ and since $\dtz\eqs \dz$, we actually have $\det_\cc A(tz)\eqs \bigl(\vol (\cp_\dz(z))\bigr)^{\frac12}$.\\
Therefore $\det_\rr (d_\Lambda \hl)\geqs\vol (\cp_\dz(z))$ and
\begin{align*}
 \iiijk&\leqs
 \int_{\over{z\in\cp_\varepsilon (z_0)\cap D}{\over{t\in[t_0,1-\gamma\dz ]}{\zeta\in \cp_\dz(z)}}} \frac1{\vol(\cp_\dz(z))} \frac{|\theta(\zeta  )|[e_j(\zeta  ),e_k(\zeta  )]}{k(\zeta  ,e_j(\zeta  ))k(\zeta  ,e_k(\zeta  ))}d\lambda(\zeta)dt d\lambda(z).
\end{align*}
Now, we proceed exactly as in the previous subsections. We integrate for $t\in[1-\gamma\dz,1-\frac\gamma2\dz]$ and get
\begin{align*}
 \iiijk&\leqs
 \int_{\over{z\in\cp_\varepsilon (z_0)\cap D}{\zeta\in \cp_\dz(z)}} \frac1{\vol(\cp_\dz(z))} \frac{\dz |\theta(\zeta  )|[e_j(\zeta  ),e_k(\zeta  )]}{k(\zeta  ,e_j(\zeta  ))k(\zeta  ,e_k(\zeta  ))}d\lambda(\zeta)d\lambda(z).
\end{align*}
We use again Fubini's theorem and get since $\dzeta\eqs\dz$ 
\begin{align*}
 \iiijk&\leqs
 \int_{\over{\zeta\in\cp_{K\varepsilon} (z_0)\cap D}{z\in K\cp_{\dzeta}(\zeta)}} \frac1{\vol(\cp_{\dzeta}(\zeta))} \frac{\dzeta|\theta(\zeta  )|[e_j(\zeta  ),e_k(\zeta  )]}{k(\zeta  ,e_j(\zeta  ))k(\zeta  ,e_k(\zeta  ))} d\lambda(z)d\lambda(\zeta).
\end{align*}
We now integrate successively for $z\in K\cp_{\dzeta}(\zeta)$ and $\zeta\in\cp_{K\varepsilon} (z_0)\cap D$ and get
\begin{align*}
\iiijk&\leqs \sigma(\cp_\varepsilon (z_0)\cap bD)\|\deltap\theta\|_{W^1_{1,1}}. 
\end{align*}
This finally ends to prove that $\iii\leqs \sigma(\cp_\varepsilon (z_0)\cap bD)\|d\hskip-2pt\cdot\hskip -2pt\theta\|_{W^1_{1,1}}$, which completes the proof of Theorem \ref{th1}.

\section{The $\overline\partial$-equation}\mlabel{section3}
The solution of the $\overline\partial$-equation will be given by the integral operator already used in \cite{DM} by K. Diederich and E. Mazzilli and which we now recall.\\
Let $\cv=\{z,\ \dz <\eta_0\}$, $\eta_0>0$, be a small neighborhood of $bD$ and let  $S\in C^\infty(\cv\times \overline D)$ be the support function constructed in \cite{DF} by K. Diederich and J. E. Forn\ae ss and globalized in \cite{WA0}. Let $Q=(Q_1,\ldots, Q_n)$ be its Hefer decomposition defined in \cite{WA} so that $S(\zeta,z)=\langle Q(\zeta,z),\zeta-z\rangle$. The support function $S$ and its Hefer decomposition are 
holomorphic in $D$ for all fixed $\zeta\in \cv$.
Let also $\chi$ be a $C^\infty$ cut-off function such that $\chi(z)=1$ if $r(z)\leq -\eta_0$ and $\chi(z)=0$ if $r(z)\geq -\frac{\eta_0}2$. We then put for $(\zeta,z)\in\overline D\times\overline D$ 
\begin{align*}
 s(\zeta,z)&:=-r(z)\sum_{i=1}^n(\overline{\zeta_i-z_i})d\zeta_i+(1-\chi(z))\overline {S(z,\zeta)}\sum_{i=1}^n Q_i(z,\zeta)d\zeta_i,\\
 q(\zeta,z)&:=\frac1{r(\zeta)}\mat({(1-\chi(\zeta  ))\sum_{i=1}^nQ_i(\zeta,z)d\zeta_i+\chi(\zeta  )\sum_{i=1}^n \diffp{r}{\zeta_i}(\zeta)d\zeta_i}),\\
 K(\zeta,z)&:=c_n\sum_{k=0}^{n-1} \frac {s(\zeta,z)\wedge (\overline\partial_\zeta s(\zeta,z))^{n-1-k}\wedge (\overline\partial_\zeta q(\zeta,z))^k}{\langle s(\zeta,z),\zeta-z\rangle^{n-k}(1-\langle q(\zeta,z),\zeta-z\rangle)^{k+1}}.
\end{align*}
We get from Berndtsson-Andersson's theorem \cite{BA}
\begin{proposition}
 Let $\omega$ be a $\overline\partial$-closed $(0,1)$-form smooth on $\overline D$. Then
\begin{align}
u(z)&:=\int_D \omega(\zeta)\wedge K(\zeta,z), \ z\in D \mlabel{u}
\end{align}
satisfies $\overline\partial u=\omega$ on $D$.
\end{proposition}
K. Diederich and E. Mazzilli showed that $K$ is uniformly integrable and get from the theorem of H. Skoda (see \cite{Sko}) that $u$ given by (\ref{u}) is continuous up to the boundary and its boundary values are still given by (\ref{u}). Following the idea of \cite{AC} also used in \cite{BG}, we prove that when $\omega$ is a smooth $\overline\partial$-closed $(0,1)$-form such that $\|\omega\|_{W_{(0,1)}^1}$ is finite, the function $\exp(pu)$, $u$ given by (\ref{u}), is in $L^1(bD)$ for some positive $p$.\\
Since $\omega(\zeta)\wedge K(\zeta,z)$ is an $(n,n)$-form, we have $\omega(\zeta)\wedge K(\zeta,z)=\psi(\zeta,z)d\lambda(\zeta)$ where $\psi(\zeta,z)=\frac1{\det(e_1,\ldots, e_n,\overline e_1,\ldots, \overline e_n)}\omega(\zeta)\wedge K(\zeta,z)(e_1,\ldots, e_n,\overline e_1,\ldots, \overline e_n)$ for any basis  $e_1,\ldots, e_n$ of $\cc^n$. In \cite{BG} they chose an $\varepsilon$-extremal basis as a basis to compute $\psi$. Here our hypothesis on $\omega$ are linked to vectors fields. That's  why we will use the same basis as in Section \ref{section2} : Let $e_j(\zeta)$ be the $j^{\rm th}$ column of the matrix $A(\zeta)=\Phi^{-1}(B(\zeta))$.

We then have
\begin{align*}
 |\psi(\zeta,z)|&\leq \sum_{i=1}^n \frac{|K(\zeta,z)(\widehat{\overline {e_i(\zeta)}})| k(\zeta,\overline {e_i(\zeta)})}{|\det(B(\zeta))|^{-1}}\frac1{k(\zeta,\overline {e_i(\zeta)})}
|\omega(\zeta)|(\overline {e_i(\zeta))})
\end{align*}
where $\widehat{\overline {e_i(\zeta)}}$ is the vectors family $e_1(\zeta),\ldots,e_n(\zeta),\overline e_1(\zeta),\ldots,\overline e_{i-1}(\zeta),\overline e_{i+1}(\zeta),$ $\ldots,$ $\overline e_{n}(\zeta)$.\\
We set $\tilde \psi_i(\zeta,z):=\frac{|K(\zeta,z)(\widehat{\overline {e_i(\zeta)}})| k(\zeta,\overline {e_i(\zeta)})}{|\det(B(\zeta))|^{-1}}$ so that
$$|u(z)|\leq \sum_{i=1}^n \int_{D} \tilde \psi_i(\zeta,z)\frac1{k(\zeta,\overline {e_i(\zeta)})}
|\omega(\zeta)|(\overline {e_i(\zeta))})d\lambda(\zeta).$$
Therefore it suffices to show that for all $i$ and all Carleson measure $\nu$, there exists $p_\nu>0$ such that for all $p<p_\nu$, the function 
$v(z):=\int_{D}\tilde\psi_i(\zeta,z)d\nu(\zeta)$, $z\in bD$, is such that $\exp( pv)$ belongs to $L^1(bD)$ and has $L^1$ norm controled by $\|\nu\|_{W^1(D)}$.\\
Now we proceed similarely to \cite{AC,BG}. We set for $f\in L^1(bD)$ 
$$L_i(f)(\zeta):=\int_{bD} \tilde\psi_i(\zeta,z)f(z)d\sigma(z),\ \zeta\in D$$
and we aim to prove the following lemma
\begin{lemma}\label{lemIII.1}
 For all $i=1,\ldots,n$
\begin{enumerate}[(i)]
 \item \mlabel{i}$\int_{bD} \tilde\psi_i(\zeta,z)d\sigma(z)\leq C$, uniformly in $\zeta\in D$,
\item \mlabel{ii}for all $\nu\in W^1(D)$, all $f\in L^1(bD)$ and all $s>0$ we have 
$$\nu\{\zeta, \ |L_i(f)(\zeta)|\geq s\}\leqs \frac{1}{s} \|\nu\|_{W^1(D)}\|f\|_{L^1(bD)}$$
 uniformly \wrt  $f$, $\nu$ and $s$.
\end{enumerate}
\end{lemma}
In order to prove Lemma \ref{lemIII.1}, we first notice as in \cite{DM}, that the denominator of $K(\zeta,z)$ is bounded away from $0$ when $\zeta$ is far from $bD$ or when $|\zeta-z|\geqs 1$, so the only case which has to be investigated is when $\zeta$ is near the boundary and it suffices to integrate for $z$ in small neighborhood of $\zeta$, say $\cp_{\varepsilon_0}(\zeta)$, $\varepsilon_0$ not depending on $\zeta$. We thus fix a point $\zeta_0$ near $bD$, $\varepsilon_0>0$ small enough, and consider points $z\in bD\cap \cp_{\varepsilon_0}(\zeta_0).$\\
We set $Q(\zeta,z):=\sum_{i=1}^n Q_i(\zeta,z)d\zeta_i$ and $\tilde Q(\zeta,z):=\sum_{i=1}^n Q_i(z,\zeta)d\zeta_i$. Let us  notice that $\tilde Q$ is holomorphic \wrt $\zeta$. So, when $z$ belong to $bD$, only the term with $k=n-1$ in the sum which defines $K$ matters. Now as in \cite{DM} we write this term as $K_1+K_2$ where 
\begin{align*}
 K_1(\zeta,z)&=\frac {\tilde Q(\zeta,z)\wedge\overline\partial r(\zeta)\wedge Q(\zeta,z)\wedge (\overline\partial_\zeta Q(\zeta,z))^{n-2}}{ S(z,\zeta)(r(\zeta)+S(\zeta,z))^{n}},\\ 
K_2(\zeta,z)&=\frac {r(\zeta)\tilde Q(\zeta,z)\wedge (\overline\partial_\zeta Q(\zeta,z))^{n-1}}{ S(z,\zeta)(r(\zeta)+S(\zeta,z))^{n}}.
\end{align*}
We will only estimate $K_1$, $K_2$ can treated with the same computations. Most of these computations are similar to those of \cite{WA1}. We write $K_1$ in a Yu-basis at $\zeta$. We first recall the definition of a Yu basis at $\zeta$ for convex domains of finite type.

The variety 1-type $\Delta_1(bD_{r(\zeta)},\zeta)$ of $bD_{r(\zeta)}$ at a point $\zeta$ is defined as
$$\Delta_1(bD,\zeta)=\sup_z\frac{\nu(z^*r)}{\nu(z-\zeta)}$$
where the supremum is taken over all non zero germ $z:\Delta \to\cc^n$ from $\Delta $, the unit disc of $\cc$, into $\cc^n$, such that $z(0)=\zeta$. The function $z^*r$ is the pullback of $r$ by $z$.\\
The variety $q$-type $\Delta _q(\zeta,bD)$ at the point $\zeta$ is then defined as
$$\Delta _q(bD,\zeta):=\inf_H \Delta _1(bD_{r(\zeta)}\cap H,\zeta)$$
where the infimum is taken over all $(n-q+1)$-dimensional complex linear manifolds $H$ passing through $\zeta$. Finally, the multitype ${\rm M}(bD_{r(\zeta)},\zeta)$ of $bD_{r(\zeta)}$ at the point $\zeta$ is defined to be the $n$-tuple $(\Delta _{n}(bD_{r(\zeta)},\zeta),\Delta _{n-1}(bD_{r(\zeta)},\zeta),\ldots, \Delta _{1}(bD_{r(\zeta)},\zeta))$.
 From Corollary 2.21 of \cite{Hef2}, we have, uniformly with respect to $\zeta$ and $\varepsilon$, $\vol(\cp_\varepsilon(\zeta))\eqs\varepsilon^{2(\Delta _{1}(bD_{r(\zeta)},\zeta)+\ldots+\Delta _{n}(bD_{r(\zeta)},\zeta))}$.

A basis $w'_1,\ldots, w'_n$ of $\cc^n$ such that for all $i$, the order of contact of $bD_{r(\zeta)}$ and the line spanned by $w'_i$ passing through $\zeta$ is equal to  $\Delta_{n+1-i}(bD_{(\zeta)},\zeta)$ is called a Yu basis at $\zeta$ (see \cite{Hef2}).

A Yu basis satisfies the following proposition which is the analog of Proposition \ref{prop0.1} for the extremal basis (see \cite{Hef2}, Theorem 2.22).
\begin{proposition}\mlabel{prop4}
Let $\zeta  \in D$ be a point near $bD$, let $(m_1,\ldots,m_n)$ denote the multitype of $bD_{r(\zeta)}$ at $\zeta$, let $w'_1,\ldots, w'_n$ be a Yu basis at $\zeta$, let $\varepsilon$ be a positive number and let $v=\sum_{j=1}^n v'_j w'_j$ be a unit vector. Then, uniformly \wrt $\zeta, v$ and $\varepsilon $ we have
\begin{align*}
\frac1{\tau (\zeta  ,v,\varepsilon )}&\eqs \sum_{j=1}^n\frac{|v'_j|}{\varepsilon ^{\frac1{m_j}}}.
\end{align*}
\end{proposition}
We notice that in particular, with the notations of Proposition \ref{prop4}, $\tau(z,w'_j,\varepsilon)\eqs \varepsilon^{\frac1{m_j}}$.

We fix a Yu basis $w'_1,\ldots,w'_n$ at $\zeta_0$ and analogously to the extremal basis notation, we put $\tau'_i(\zeta,\varepsilon):=\tau(\zeta,w'_i,\varepsilon )$. We denote by $\zeta'=(\zeta'_1,\ldots,\zeta'_n)$ the coordinates of a point $\zeta$ in the coordinates system centered at $\zeta_0$ of basis $w'_1,\ldots,w'_n$.
Then we write $\tilde Q$ and $Q$ in the Yu basi at $\zeta_0$:
$\tilde Q(\zeta,z)=\sum_{j=1}^n\tilde Q'_j(\zeta,z)d\zeta'_j$ and 
$Q(\zeta,z)=\sum_{j=1}^nQ'_j(\zeta,z)d\zeta'_j$. 
In the Yu basis $w'_1,\ldots, w'_n$, $K_1$ is a sum of the following terms
\begin{align*}
{K_{\nu,\mu}(\zeta,z)}={\frac {\tilde Q'_{\nu_1}(\zeta,z)d\zeta'_{\nu_1}\wedge
\diffp{ r}{\overline \zeta'_{\mu_2}}(\zeta)d{\overline\zeta'_{\mu_2}}\wedge
 Q'_{\nu_2}(\zeta,z)d\zeta'_{\nu_2}
 \bigwedge_{i=3}^{n}\diffp {Q'_{\nu_i}}{\overline\zeta'_{\mu_i}}(\zeta,z)d\overline\zeta'_{\mu_i}\wedge d\zeta_{\nu_i}'}{ S(z,\zeta)(r(\zeta)+S(\zeta,z))^{n}} }&
\end{align*}
where $\nu_i$ and $\mu_i$ run from 1 to $n$, $\nu_i\neq \nu_j$, $\mu_i\neq\mu_j$ for $i\neq j$. We have to estimate 
 $\frac{k(\zeta,\overline {e_{i_0}(\zeta)})}{|\det(B(\zeta))|^{-1}}K_{\nu,\mu}(\zeta,z)(\widehat{\overline{e_{i_0}(\zeta)}})$ for all such $\nu$ and $\mu$.
We have the following proposition which comes from \cite{WA,WA3} and Proposition \ref{prop4}, and which were already used in \cite{WA1} :
\begin{proposition}\mlabel{prop6}
For all $\zeta$ near enough $bD$, all sufficiently small $\varepsilon>0$, all $\xi,z\in\cp_\varepsilon(\zeta)$ and $i, j =1,\ldots, n$, we have uniformly with respect to $\zeta$, $z$, $\xi$  and $\varepsilon$ 
\begin{align*}
|Q'_i(\xi,z)|&\leqs \frac{\varepsilon}{\tau'_j(\zeta,\varepsilon)},\\
 \mat|{
\diffp{Q'_i} {\overline {\zeta'_j}} (\xi,z) }|&\leqs \frac{\varepsilon} {\tau'_j(\zeta,\varepsilon)\tau'_i(\zeta,\varepsilon)}.
\end{align*}
\end{proposition}

Since $\delta$ is a pseudodistance, we deduce the following inequalities from Proposition 4.4 from \cite{WA1}  : for all $z\in \cp_{\varepsilon}(\zeta_0)\setminus c\cp_\varepsilon(\zeta_0)\cap bD$, $c>0$ given by Corollary \ref{corI.1.2} such that $c\cp_\varepsilon(\zeta_0)\subset \cp_{\frac12\varepsilon}(\zeta_0)$, we have uniformly \wrt $\zeta_0$ and $z$
\begin{align}
|S(z,\zeta_0)|&\geqs\varepsilon,\mlabel{eq7}\\
|S(\zeta_0,z)+r(\zeta_0)|&\geqs\varepsilon\mlabel{eq8}.
\end{align}

With the inequality $\mat|{\diffp{ r}{\overline\zeta'_{\mu_2}}(\zeta_0)}|\leqs\frac{\dzetaO }{\tau'_{\mu_2}(\zeta_0,\dzetaO )}$ which comes from proposition \ref{propI.1.8}, we get

\begin{align}
&\nonumber {\mat|{\frac {\tilde Q'_{\nu_1}(\zeta_0,z)\diffp{ r}{\overline\zeta'_{\mu_2}}(\zeta_0)
Q'_{\nu_2}(\zeta_0,z)
 \prod_{i=3}^{n}\diffp {Q'_{\nu_i}}{\overline\zeta'_{\mu_i}}(\zeta_0,z)}{ S(z,\zeta_0)(r(\zeta_0)+S(\zeta_0,z))^{n}} }|}\\&\hskip200pt \leqs{\frac{\dzetaO }{\tau'_{\mu_2}(\zeta_0,\dzetaO )}\frac1{\varepsilon\prod_{i=1}^{n}\hskip -3pt \tau'_{i}(\zeta_0,\varepsilon\hskip -1pt)\prod_{i=3}^{n}\tau'_{\mu_i}(\zeta_0,\varepsilon)}}\mlabel{eq2}
\end{align}
We now estimate $\frac{k(\zeta_0,\overline{e_{i_0}(\zeta_0)})}{(\det B(\zeta_0))^{-1}}\bigwedge_{i=1}^n d\zeta'_i\wedge \bigwedge_{i=2}^n d\overline\zeta'_{\mu_i}(\widehat{\overline e_{i_0}(\zeta_0)})$.
We denote the coordinates of $e'_i(\zeta_0)$ in the Yu basis $w'_1,\ldots, w'_n$ by $(e'_{i1}(\zeta_0),\ldots,e'_{in}(\zeta_0))$. Proposition 4.8 of \cite{WA1} asserts that
$$|e'_{ij} (\zeta_0)|\leqs \tau'_j(\zeta_0,\dzetaO).$$
Therefore 
\begin{align*}
\mat|{\bigwedge_{i=1}^n d\zeta'_i\wedge \bigwedge_{i=2}^n d\overline\zeta'_{\mu_i} (\widehat{\overline e_{i_0}(\zeta_0)})}|\leqs\prod_{i=1}^n\tau'_i(\zeta_0,\dzetaO) \prod_{i=2}^n\tau'_{\mu_i}(\zeta_0,\dzetaO).
\end{align*}
Since $e_1(\zeta_0),\ldots,e_n(\zeta_0)$ is an orthonormal basis for the Bergman metric, $k(\zeta_0,\overline{e_{i_0}(\zeta_0)})\eqs\dzetaO$ and since $(\det B(\zeta_0))^{-1} \eqs\vol(\cp_{\dzetaO}(\zeta_0))\eqs\prod_{i=1}^n\tau'_i(\zeta_0,\dzetaO)^2$, we get
\begin{align}
 \left|\frac{k(\zeta_0,\overline{e_{i_0}(\zeta_0)})}{(\det B(\zeta_0))^{-1}}\bigwedge_{i=1}^n d\zeta'_i\wedge \bigwedge_{i=2}^n d\overline\zeta'_{\mu_i}(\widehat{\overline{e_{i_0}(\zeta)}})\right|
 &\leqs \frac{\dzetaO}{\tau'_{\mu_1}(\zeta_0,\dzetaO)},\mlabel{eq1}
\end{align}
where $\mu_1\in\{1,\ldots,n\}\setminus\{\mu_2,\ldots, \mu_n\}$.\\
From  (\ref{eq2}) and  (\ref{eq1}), we obtain
\begin{align*}
\mat|{\frac{k(\zeta_0,\overline {e_{i_0}(\zeta_0)})}{\det(B(\zeta_0)^{-1})}K_{\nu,\mu}(\zeta_0,z)(\widehat{\overline e_{i_0}(\zeta_0)}) }|\leqs 
\frac{\dzetaO^2 \tau'_{\mu_1}(\zeta_0,\varepsilon)\tau'_{\mu_2}(\zeta_0,\varepsilon)}{\varepsilon^2 \tau'_{\mu_1}(\zeta_0,\dzetaO)\tau'_{\mu_2}(\zeta_0,\dzetaO)}\frac1{{\sigma(\cp_\varepsilon(\zeta_0)\cap bD)}}.
\end{align*}
Since at least $\mu_1$ or $\mu_2$ is different from 1, we get from Proposition \ref{prop4}
\begin{align}
\frac{\dzetaO^2 \tau'_{\mu_1}(\zeta_0,\varepsilon)\tau'_{\mu_2}(\zeta_0,\varepsilon)}{\varepsilon^2 \tau'_{\mu_1}(\zeta_0,\dzetaO)\tau'_{\mu_2}(\zeta_0,\dzetaO)}&\leqs 
\mat({\frac{\dzetaO }{\varepsilon}})^{\frac12}.
\end{align}
Finally, we get for all $z\in bD\cap \cp_{\varepsilon}(\zeta_0)\setminus c\cp_\varepsilon(\zeta_0)$
\begin{align}
\mat|{\frac{k(\zeta_0,\overline {e_{i_0}(\zeta_0)})}{\det(B(\zeta_0)^{-1})}K_1(\zeta_0,z)(\widehat{\overline e_{i_0}(\zeta_0)}) }|\leqs 
\mat({\frac{\dzetaO }{\varepsilon}})^{\frac12}\frac1{\sigma(\cp_\varepsilon(\zeta_0)\cap bD)}.\mlabel{eq10}
\end{align}

Since $z$ belongs to $bD$ and since $\zeta_0$ belongs to $D$, $\delta(z,\zeta_0)\geqs d(\zeta_0)$ so the inequalities (\ref{eq7}) and (\ref{eq8}) are still valid for $z\in \cp_{\dzetaO }\cap bD$. Therefore we have for such $z$ :
\begin{align}
\mat|{\frac{k(\zeta_0,\overline {e_{i_0}(\zeta_0)})}{\det(B(\zeta_0)^{-1})}K_1(\zeta_0,z)(\widehat{\overline e_{i_0}(\zeta_0)}) }|\leqs 
\frac1{\sigma(\cp_{\dzetaO }(\zeta_0)\cap bD)}.\mlabel{eq9}
\end{align}
The estimates (\ref{eq10}) and (\ref{eq9}) can be shown for $K_2$ instead of $K_1$.
Now, as in \cite{DFF} we cover $\cp_{\varepsilon_0}(\zeta_0)$ 
 with some polyannuli based on McNeal's polydiscs. For sufficiently small $\varepsilon>0$ we set
$\cp^i_\varepsilon(\zeta_0):= \cp_{2^{i}\varepsilon}(z_0)\setminus c \cp_{2^{i}\varepsilon}(z_0).$
This gives us the following covering 
\begin{eqnarray}
\cp_{\varepsilon_0}(\zeta_0)\subset \cp_{\dzetaO }(\zeta_0) \cup\bigcup_{k=0}^{k_0} \cp^k_{\dzetaO }(\zeta_0) \mlabel{covering}
\end{eqnarray}
where $k_0$  satisfies $k_0\eqs |\ln\varepsilon_0- \ln \dzetaO |$, uniformly in $\zeta_0$ 
 and $\varepsilon_0$. We finally get
\begin{align*}
 \lefteqn{\int_{bD\cap \cp_{\varepsilon_0}(\zeta_0)} \tilde\psi_{i_0}(\zeta_0,z)d\sigma(z)\leqs}\\
&\leqs
\int_{bD\cap \cp_{\dzetaO }(\zeta_0)}\mat|{\frac{k(\zeta_0,\overline {e_{i_0}(\zeta_0)})}{\det(B(\zeta_0)^{-1})}K(\zeta_0,z)(\widehat{\overline e_{i_0}(\zeta_0)}) }| d\sigma(z)\\
&+\sum_{k=0}^{k_0}\int_{bD\cap \cp^k_{\dzetaO }(\zeta_0)}\mat|{\frac{k(\zeta_0,\overline {e_{i_0}(\zeta_0)})}{\det(B(\zeta_0)^{-1})}K(\zeta_0,z)(\widehat{\overline e_{i_0}(\zeta_0)}) }| d\sigma(z)\\
&\leqs 1+\sum_{k=0}^{k_0}\frac1 {2^{k}}\leqs 1
\end{align*}
uniformly \wrt $\zeta_0$. The first point of Lemma \ref{lemIII.1} is then proved.\\
In order to prove the second point, for $f\in L^1(bD)$ we defined the maximal function
$$\tilde f(\zeta):=\sup_{\varepsilon\geq \dzeta}\left( \frac1{\sigma(\cp_\varepsilon(\zeta)\cap bD)}
\int_{\cp_\varepsilon(\zeta)\cap bD}|f(\xi)|d\sigma(\xi)\right).$$
We deduce from (\ref{eq10}) and (\ref{eq9}) and their analogous version for $K_2$ that
\begin{align*}
 \lefteqn{\int_{bD\cap \cp_{\varepsilon_0}(\zeta_0)}\mat|{\frac{k(\zeta_0,\overline {e_{i_0}(\zeta_0)})}{\det(B(\zeta_0)^{-1})}K(\zeta_0,z)(\widehat{\overline e_{i_0}(\zeta_0)}) }||f(z)| d\sigma(z)\leqs}\\ 
&\leqs
\int_{bD\cap \cp_{\dzetaO }(\zeta_0)}\mat|{\frac{k(\zeta_0,\overline {e_{i_0}(\zeta_0)})}{\det(B(\zeta_0)^{-1})}K(\zeta_0,z)(\widehat{\overline e_{i_0}(\zeta_0)}) }||f(z)| d\sigma(z)\\
&+\sum_{k=0}^{k_0}\int_{bD\cap \cp^k_{\dzetaO }(\zeta_0)}\mat|{\frac{k(\zeta_0,\overline {e_{i_0}(\zeta_0)})}{\det(B(\zeta_0)^{-1})}K(\zeta_0,z)(\widehat{\overline e_{i_0}(\zeta_0)}) }||f(z)| d\sigma(z)\\
&\leqs \frac1{\sigma(\cp_{\dzetaO }(\zeta_0)\cap bD)}\int_{bD\cap \cp_{\dzetaO }(\zeta_0)}|f(z)|d\sigma(z)\\
&+\sum_{k=0}^{k_0}\frac1{2^k}\frac1{\sigma(\cp_{2^k\dzetaO }(\zeta_0)\cap bD)}\int_{bD\cap \cp_{2^k\dzetaO }(\zeta_0)}|f(z)|d\sigma(z)\\
&\leqs \tilde f(\zeta_0).
\end{align*}
Now Theorem 2.1 of \cite{Hor} gives 
$$\nu\{\zeta, |L_i(f)(\zeta)| > s\}\leq \nu\{\zeta, |\tilde f(\zeta)| > \frac s C\}\leqs \frac{1}{s} \|\nu\|_{W^1(D)}\|f\|_{L^1(bD)},$$ uniformly \wrt $f$ and $\nu$. This ends the proof of Lemma \ref{lemIII.1}.

We now prove that $v(z):=\int_{D}\tilde\psi_i(\zeta,z)d\nu(\zeta)$ is such that $\exp (pv)$ belongs to $L^1(bD)$ for some $p>0$ depending only on $\|\nu\|_{W^1}$. The method is exactly the same as in \cite{AC,BG}. We include it for completeness.\\
Let $E_t:=\{z\in bD,\ v(z)>t\}$. We have
$$\int_{bD} \exp(pv(z))d\sigma(z)=\int_0^\infty p\exp({pt})\sigma(E_t)dt+\sigma(bD).$$
We claim that there exist $C,C'>0$ not depending on $t$ or $\nu$ such that $\sigma(E_t)\leq C  e^{-\frac {C't}{\|\nu\|_{W^1}}}$.
Then, for $p<\frac {C'}{{\|\nu\|_{W^1}}}$, $\int_{bD} \exp(pv(z))d\sigma(z)$ is bounded, which was to be shown. Theorem \ref{th2} will therefore be proved as soon as the claim is proved.\\
We denote by $\chi_{E_t}$ the characteristic function of $E_t$. We have 
\begin{align*}
 t\sigma(E_t)&\leq \int_{E_t} v(z)d\sigma(z)\\
&\leq \int_{D}\mat({\int_{bD}\tilde\psi_i(\zeta,z)\chi_{E_t}(z)d\sigma(z)})d\nu(\zeta)\\
&\leq \int_DL_i(\chi_{E_t})(\zeta)d\nu(\zeta)\\
&\leq \int_0^\infty \nu(\{L_i(\chi_{E_t})>s\})ds.
\end{align*}
Lemma \ref{lemIII.1} (\ref{i}) implies that $L_i(\chi_{E_t})$ is bounded by some constant $M>0$ which does not depend on $t$ so
\begin{align*}
t\sigma(E_t)&\leq \int_0^{M} \nu(\{L_i(\chi_{E_t})>s\})ds\\
&\leq \int_0^{\sigma(E_t)}\nu(\{ L_i(\chi_{E_t})>s\})ds+\int_{\sigma(E_t)}^{M}\nu(\{ L_i(\chi_{E_t})>s\})ds.
\end{align*}
Now Lemma \ref{lemIII.1} (\ref{ii}) yields
\begin{align*}
t\sigma(E_t)&\leqs \sigma(E_t)\sigma(bD)\|\nu\|_{W^1}+\int_{\sigma(E_t)}^{M}\|\nu\|_{W_1}\frac 1s \sigma(E_t)ds \\
&\leqs \sigma(E_t)\|\nu\|_{W^1}\left(\sigma(bD)+\ln \mat({\frac{M}{\sigma(E_t)} })\right).
\end{align*}
Therefore there exists $C',C>0$ which does not depend on $\nu$ or $t$ such that $\sigma(E_t)\leq Ce^{-\frac {C't}{\|\nu\|_{W^1}}}$ and the claim is proved.

\end{document}